%% file: logahoric.tex
\documentclass [12pt]{amsart}

\usepackage{amssymb,amsxtra,amsfonts}
\usepackage{graphics}

\usepackage[colorlinks]{hyperref}

\input{macros}

\input{macros-thm1}
\renewenvironment{ppb}[1]{}{}

\begin{document}

\title{Riemann--Hilbert for tame complex parahoric connections}
\author{Philip Boalch}%

\begin{abstract}
A local Riemann--Hilbert correspondence for tame meromorphic connections on a curve compatible with a parahoric level structure will be established.
Special cases include logarithmic connections on $G$-bundles and 
on parabolic $G$-bundles.
The corresponding Betti data involves pairs $(M,P)$ consisting of the local monodromy $M\in G$ and a (weighted) parabolic subgroup $P\subset G$ such that $M\in P$, as in the multiplicative Brieskorn--Grothendieck--Springer resolution (extended to the parabolic case).
We will also construct the natural quasi-Hamiltonian structures that arise on such spaces of enriched monodromy data.
\end{abstract}

\maketitle

\section{Introduction}

The starting point of this article was an attempt to extend to $G$-bundles the local classification of logarithmic connections on vector bundles on curves in terms of Levelt filtrations, where $G$ is a connected complex reductive group.
Namely logarithmic connections on vector bundles are classified locally
by triples $(V,F,M)$ where
$V$ is a finite dimensional complex vector space,
$F$ is a decreasing finite filtration of $V$ indexed by $\IZ$
and $M\in \GL(V)$ preserves the filtration $F$.
\ppb{
\begin{align*}
\bullet & \text{ $V$ is a finite dimensional complex vector space}\\
\bullet & \text{ $F$ is a decreasing finite filtration of $V$ indexed by $\IZ$}\\
\bullet & \text{ $M\in \GL(V)$ preserves the filtration $F$.}\\
\end{align*}
}
If we forget the filtration then we obtain the local classification of regular singular connections, much studied e.g. by Deligne \cite{Del70} (in arbitrary dimensions)---they form a Tannakian category (cf. \cite{katz-dgps}) and the extension to $G$-bundles is then straightforward
(they are classified by their monodromy $M\in G$ up to conjugation) although a direct approach is possible (see \cite{BV83}).

Thus for general $G$ we wish to describe the extra data needed to determine a logarithmic connection and establish the precise correspondence.
Unfortunately the category of triples $(V,F,M)$ is not abelian, and so not Tannakian, and so it seems a direct approach is necessary (if it were Tannakian we could just take the space of homomorphisms from the corresponding group into $G$).
The key point in the above classification of logarithmic connections %
is that one  may choose a local holomorphic trivialisation and a one-parameter subgroup $\varphi:\IC^*\to \GL(V)$ such that if we view $\varphi$ as a meromorphic gauge transformation, then in the resulting trivialisation the connection  takes the simple form
$$R\frac{dz}{z}$$
 for some $R\in \End(V)$ with eigenvalues all 
having real parts in the interval  $[0,1)$ (using a fixed local coordinate $z$).
The resulting data is then $(V,F,M)$ where $M=e^{2\pi i R}$ is the local monodromy and $F$ is the filtration naturally associated to $\varphi$.
The utility of the filtration is that if $g\in \GL(V)$ and $\psi=g\varphi g^{-1}$ is a conjugate one parameter subgroup then the meromorphic group element $\varphi\psi^{-1}$ is holomorphic if and only if $\varphi$ and $\psi$ determine the {\em same} filtration, i.e. $g$ preserves $F$.
This is why the Levelt filtration (from \cite{levelt61} (2.2)) gives a much cleaner approach than the naive viewpoint of directly recording the extra terms that may occur in the case of ``resonant'' connections.

For general $G$ the notion of flag generalises directly to the notion of parabolic subgroup, and one may in general attach a parabolic subgroup $P\subset G$ to a one parameter subgroup (see e.g. Mumford \cite{MFK94} p.55). 
However it is {\em not true}, even for $\SL_2(\IC)$, that every logarithmic connection may be put in the simple form $Rdz/z$ with $R\in \g=\Lie(G)$ via a suitable trivialisation and a one parameter subgroup (see \cite{BV83} p.65),
and even if we did restrict to such connections a good analogue of the above normalisation of the eigenvalues looks to be  elusive.
At first sight this is bad news since it means the direct analogue of the above $\GL_n(\IC)$ classification does not seem to hold, but it is also good news: the failure to reduce to the simple form corresponds directly to the fact that there are logarithmic connections whose monodromy $M$ is not in the image of the exponential map, so we can hope for a more complete correspondence involving all possible monodromy conjugacy classes.

Whilst extending the nonabelian Hodge correspondence to open curves Simpson \cite{Sim-hboncc} gave an alternative approach, which he also applies to more general objects (``filtered tame $\cD$-modules''), but still in the context of vector bundles. In the case of logarithmic connections this amounts to refining the Levelt filtration to take into account the exact rate of growth of solutions rather than its integer part as was effectively done above.
It is this approach that we are able to extend to all complex reductive groups.
Moreover the final version of the correspondence (Theorem D) involves some new features which do not occur in the case of vector bundles.
Also a surprisingly clean statement (Corollary E) is possible if we use Bruhat--Tits buildings.

Our  motivation was to understand the  spaces of monodromy type data  that occur in the extension  of 
the nonabelian Hodge correspondence to the case of irregular connections on curves \cite{Sab99, wnabh}, and its extension to arbitrary $G$. Using the quasi-Hamiltonian approach this problem may be broken up into pieces: understanding the Stokes data, and understanding what to do for regular singularities. Since it is possible to understand the Stokes data for arbitrary $G$ (cf. \cite{bafi, saqh, qf}) we are left with the problem of extending Simpson's tame Riemann--Hilbert correspondence \cite{Sim-hboncc} to general $G$, which we will do here.
At the end of the day this will give the algebraic ``Betti'' description of some complex manifolds supporting hyperk\"ahler metrics
(appearing in the nonabelian Hodge theory of  curves).
Some motivation also came from trying to understand
 the recent work of Gukov--Witten \cite{gukov-witten, gukov-witten-rigid} on the tamely ramified geometric Langlands correspondence
(in particular this justifies our desire to work uniformly with arbitrary complex reductive groups).

\newpage
\ssn{Results and further evolution}

We will state three local classification results, of increasing complexity, since each may be of interest to different readers.
In the case of logarithmic connections the statement is as follows.
Let $\lt\subset \g$ be a Cartan subalgebra corresponding to a maximal torus $T\subset G$, and let $\lt_\IR = X_*(T)\otimes_\IZ \IR$ be the space of real cocharacters so that $\lt = \lt_\IR\otimes_\IR \IC$.
Choose an element $\tau+\si\in \lt$ with real part $\tau\in \lt_\IR$ and a nilpotent element $n\in \g$ commuting with $\tau+\si$.
Let $\orbit\subset \g$ be the adjoint orbit  of $\tau + \si + n \in \g.$
Let $L\subset G$ be the centraliser of $\tau$ and 
let $P_\tau\subset G$ 
be the parabolic subgroup determined by $\tau$ (see Section \ref{sn: defns}), so that $L$ is a Levi subgroup of $P_\tau$.
Let  $\cC\subset L$ be the conjugacy class containing the element
$\exp(2\pi i (\tau+\si+n))\in L.$
Then $\cC$ canonically determines a conjugacy class in the Levi factor of any parabolic subgroup of $G$ conjugate to $P_\tau$ (see Lemma \ref{lem: conjorbits}).

\vspace{0.3cm}

\noindent{\bf Theorem A} (Logarithmic case).
{\em 
There is a canonical bijection between isomorphism classes of germs of logarithmic connections on $G$-bundles with residue in $\orbit$ and conjugacy classes of pairs $(M,P)$ with 
$P\subset G$ a parabolic subgroup conjugate to $P_\tau$ and $M\in P$
such that $\pi(M)\in \cC$, where $\pi$ is the natural projection from $P$ onto its Levi factor.
}

\vspace{0.3cm}

Note that if $\orbit$ is nonresonant (i.e. $\al(\tau+\si)$ is not a nonzero integer for any root $\al$) 
then the condition  $\pi(M)\in \cC$ implies that $M$ itself is conjugate to $\exp(2\pi i (\tau+\si+n))$.

At this point we investigated the spaces of enriched monodromy data that start to appear here from a quasi-Hamiltonian viewpoint.
In the case of compact groups, when studying moduli space of flat connections on open Riemann surfaces, one fixes the conjugacy class of monodromy around each boundary component/puncture in order to obtain symplectic moduli spaces. As in \cite{Sim-hboncc} for $\GL_n(\IC)$ we now see this is not the most general thing that arises in the case of complex reductive groups:  in general one should fix the conjugacy class 
of the image in a Levi factor. In the quasi-Hamiltonian approach where one constructs spaces of (generalised) monodromy data by fusing together some basic pieces this corresponds to a ``new piece'', as follows.

Let 
$P_0\subset G$ be a fixed parabolic subgroup with Levi factor $L$.
Choose a conjugacy class $\cC\subset L$ (as remarked above this canonically determines a conjugacy class in the Levi factor of any conjugate parabolic subgroup).
Let $\IP\cong G/P_0$ be the set of parabolic subgroups conjugate to $P_0$.

\vspace{0.3cm}

\noindent{\bf Theorem B}.
{\em 
The smooth variety $\wh \cC$  of pairs $(M,P)\in G\times \IP$ such that $M\in P$ and $\pi(M)\in \cC$ is a quasi-Hamiltonian $G$-space with moment map given by $$(M,P)\mapsto M\in G.$$  
}

\vspace{0.3cm}

If $P_0$ is a Borel, these spaces appear in the multiplicative 
Brieskorn--Grothendieck--Springer resolution. If $P_0=G$ then $\wh \cC=\cC$. 
The additive analogue (on the Lie algebra level) of this is well-known, when the resolution is the moment map in the usual sense
(see \cite{biquard-nahm} Theorem 2). Some Poisson aspects of the multiplicative case are studied in \cite{evenslu-groth}, but the quasi-Hamiltonian (or quasi--Poisson) viewpoint looks to be more natural. 
The $\GL_n(\IC)$ case may be constructed differently via quivers (cf. \cite{yamakawa-mpa}).

This enables us to construct lots of complex symplectic manifolds of `enriched monodromy data' of the form
$$\left(\ID\fus\cdots\fus \ID \fus \wh \cC_1 \fus \cdots \fus \wh \cC_m\right)\spq G$$
where $\ID\cong G\times G $ is the internally fused double, the $\cC_i$ are conjugacy classes in Levi factors of various parabolic subgroups of  $G$ and ``$\spq$'' denotes a quasi-Hamiltonian quotient (a quotient of a subvariety). The problem now is to try to interpret these spaces as spaces of meromorphic connections on Riemann surfaces (of genus equal to the number of factors of $\ID$ appearing here). This almost immediately reduces to the local problem of interpreting the spaces  $\wh \cC$---clearly only some of them arise in Theorem A since  $\tau$ determines the parabolic subgroup $P_\tau$ and also arises in the choice of $\cC$. 

The next generalisation is to consider logarithmic connections on parabolic bundles as follows.
We will say an element  $\th\in \lt_\IR$ is {\em small} if $\al(\th)< 1$ for all roots $\al$.
Choose a small element $\th$ and let $P_\th\subset G$ be the corresponding 
parabolic subgroup.
A (germ of a) parabolic bundle with weight $\th$  is a $G$-bundle $E$ on a disc together with a reduction of structure group\footnote{this is a choice of a point of $E_0/P_\th$ where $E_0\cong G$ is the fibre of $E$ at $0$. Equivalently it is the choice of a parabolic subgroup conjugate to $P_\th$ in $G(E)_0\cong G$, where $G(E)$ is the associated adjoint group bundle.}
 to $P_\th$ at $0$.
A logarithmic connection on a parabolic $G$-bundle $E$ is then a logarithmic connection whose residue preserves the parabolic structure.
In local coordinates and trivialisation this means the connection takes the form
$$A= \left(\sum_{i\ge 0}A_i z^i \right)\frac{dz}{z}$$
with $A_i\in \g$ and the reduction determines a parabolic subalgebra $\lp\subset \g$ and the compatibility condition means $A_0\in \lp$.
 
The parabolic correspondence is then as follows.
Fix $\tau+\si+n\in \g$ as above and suppose further that $n$ commutes with $\th$. 
Let $H_\th\subset G$ be the centraliser of $\th$ (a Levi subgroup of $P_\th$) and let $\orbit\subset \h_\th$ be the adjoint orbit of
$\tau+\si+n\in \h_\th:=\Lie(H_\th)$.
The orbit $\orbit$ canonically determines an adjoint orbit in the Levi factor $\h$ of any parabolic subalgebra $\lp$ conjugate to $\Lie(P_\th)$. We will say a parabolic connection ``lies over $\orbit$'' if its residue (in $\lp$) projects to an element of $\orbit\subset \lh$
under the canonical map $\lp\twoheadrightarrow \lh$, quotienting by the nilradical.
Now set 
$$\phi = \tau + \th \in \lt_\IR$$
and let $P_\phi\subset G$ be the corresponding parabolic subgroup and $L\subset P_\phi$ be the centraliser in $G$ of $\phi$ (a  Levi subgroup of $P_\phi$). Then $\exp(2\pi i (\tau+\si +n))$ is in $L$ and we let $\cC\subset L$ be its conjugacy class.

\vspace{0.3cm}

\noindent{\bf Theorem C} (Parabolic case).
{\em 
Suppose that the centraliser in $G$ of $\exp(2\pi i \th)\in G$ is connected.
Then there is a canonical bijection between isomorphism classes of germs of parabolic  connections on $G$-bundles  with weight $\th$ and residue lying over $\orbit$, and conjugacy classes of pairs $(M,P)$ with 
$P\subset G$ a parabolic subgroup conjugate to $P_\phi$ and $M\in P$
such that $\pi(M)\in \cC$, where $\pi$ is the natural projection from $P$ onto its Levi factor.
}

\vspace{0.3cm}

This clearly captures many more of the spaces $\wh \cC$, and specialises to Theorem A if $\th=0$.
But it is still not entirely satisfactory for several reasons. 
First, by definition $\cC\subset L$ is always in the image of the exponential map (so we do not always get all possible classes).
Secondly Theorem C involves a connected centraliser condition---this holds automatically if the derived subgroup of $G$ is simply-connected (e.g. for $\GL_n(\IC)$ or for any simply-connected semisimple group), 
but not always. For  example Theorem C does not apply to $\PGL_2(\IC)$ and $\th=\bsmx 1 & \\ & 0 \esmx/2.$
Thirdly we have restricted to {\em small} weights $\th$ (such that $\al(\th)< 1$ for all roots $\al$).\footnote{Note for $\GL_n(\IC)$ one can always reduce to the case of small weights, every Levi subgroup has surjective exponential map, and the centraliser of any semisimple group element is connected.}

Somewhat miraculously all the problems disappear if we  pass to the objects which naturally appear when we do not restrict to small weights and if we use their most natural groups of automorphisms.
This is most simply described in local coordinates/trivialisations.
Given any $\th\in \lt_\IR$ we have a decomposition 
$\g=\bigoplus \g_\la$ of the Lie algebra of $G$ into the eigenspaces of $\ad_\th$
and we may consider the space of ``tame parahoric'' connections of the form
$$\cA_\th = \left\{ A=\left(\sum_{i\in \IZ, \la\in \IR} A_{i \la} z^i\right)\frac{dz}{z} \ \Bigl\vert\ 
A_{i \la} \in \g_\la\text{ and } i+\la\ge 0\right\}\subset 
\g(\!(z)\!)dz.$$
This is acted on (by gauge transformations) by the extended parahoric subgroup
$$\wh \cP_\th = \{ g\in G(\!(z)\!) \st z^\th g z^{-\th} \text{ has a limit as $z\to 0$ along any ray}\}$$
where $z^\th=\exp(\th\log(z))$ (see section \ref{sn: defns}).
\ppb{
Thus if $\th$ is small  $\cA_\th$ is  the space of parabolic connections with weight $\th$, and if moreover the centraliser of $\exp(2\pi i \th)$ is connected then $\wh \cP_\th$ is just the group of holomorphic bundle automorphism preserving the parabolic structure (i.e. holomorphic and with value at $z=0$ in $P_\th$).
}
The main result is the classification of $\wh \cP_\th$ orbits of such connections.
That this is a nontrivial generalisation is clear if we consider for example the case $G=E_8$: then there are $511$ conjugacy classes of parahoric subgroups, of which only $256$ arise in the parabolic case.
To describe the classification we will first discuss the generalisation of the notion of fixing the adjoint orbit of the residue.

Let $\wh H_\th\subset G$ be the centraliser of $\exp(2\pi i \th)$
(which might be disconnected),  and now set $\lh_\th = \Lie(\wh H_\th)$, which agrees with the previous definition for small $\th$.
The group $\wh H_\th$ is isomorphic to the ``Levi'' subgroup $\wh \cL_\th = \{ z^{-\th}h z^\th \st h\in \wh H_\th\}$ of $\wh \cP_\th$.
The finite dimensional weight zero piece $\cA_\th(0) = \{ \sum A_{i, -i} z^i dz/z\}$ of $\cA_\th$ is acted on by $\wh \cL_\th$ and the orbits correspond to adjoint orbits of $\wh H_\th$ (see Lemma \ref{lem: local orbits}).
The generalisation of fixing the adjoint orbit of the 
(Levi quotient of the) residue is to fix the adjoint orbit $\orbit \subset \lh_\th$ corresponding to  the weight zero part of the connection.
Notice that in general one now gets a richer class of subalgebras 
$\lh_\th\subset \g$: it is not necessarily the Levi factor of a parabolic (e.g. 
if $G=G_2$ one may obtain $\sl_3(\IC)\subset \g$ which is still simple of rank two). 
The full statement of the local correspondence is then as follows.

Fix elements $\th, \tau \in \lt_\IR$ and $\si\in \sqrt{-1}\lt_\IR$ and set $\phi=\th+\tau$.
Choose a nilpotent element $n\in \lh_\th\subset \g$ 
commuting with $\phi$ and $\si$.
(Thus there is a finite decomposition $n=\sum a_i$ 
with $[\tau,a_i]=ia_i=[a_i,\th]$ for $i\in \IZ$.)
Let $\orbit\subset \lh_\th$ be the adjoint orbit %
of the element $\phi + \si + n \in \lh_\th.$
This corresponds to the element $(\tau + \si + \sum a_i z^i)dz/z \in \cA_\th(0).$ 
Let $L\subset P_\phi$ be the Levi subgroup as above, but define $\cC\subset L$ to be the conjugacy class containing the element
$$\exp(2\pi i (\tau+\si))\exp(2\pi i n )\in L.$$
Note that $\cC$ is not necessarily an exponential conjugacy class, since $n$ and $\tau$ might not commute---indeed the Jordan decomposition implies all conjugacy classes arise in this way.

\vspace{0.2cm}

\noindent{\bf Theorem D} (Parahoric case).
{\em 
There is a canonical bijection between the $\wh \cP_\th$ orbits of 
tame parahoric connections in $\cA_\th$ lying over $\orbit$ and conjugacy classes of pairs $(M,P)$ with $P\subset G$ a parabolic subgroup conjugate to  $P_\phi$ and $M\in P$ such that $\pi(M)\in\cC$. 
}

\vspace{0.2cm}

This is the main result and specialises to Theorems A and C.
Finally, by considering the space $\IB(G)$ of weighted parabolic subgroups of $G$, and  the space $\cB(LG)$ of weighted parahoric subgroups of the local loop group $LG=G(\!( z)\!)$, it is possible to deduce the following statement, not involving orbit choices etc:

\vspace{0.2cm}

\noindent{\bf Corollary E}.
{\em 
There is a canonical bijection between $LG$ orbits of tame parahoric connections and $G$ orbits of enriched monodromy data:
$$\Bigl\{ (A,p) \st p\in \cB(LG), A\in \cA_p\Bigr\}/LG \quad \cong\quad 
 \Bigl\{ (M,b) \st b\in \IB(G), M\in P_b\Bigr\}/G.$$
}
\ppb{
For the global result, 
let $\G\to \Si$ be a parahoric (Bruhat--Tits) group scheme over a smooth complex  algebraic curve $\Si$ such that locally 
$\G$ looks like a parahoric subgroup $\cP$ of the local loop group 
$LG$ and at all but finitely many points $D\subset \Si$ it looks like $G(\cO)\cong G\flb z\frb$.
Then we consider {\em parahoric bundles} over $\Si$ i.e. 
$\G$-torsors over $\Si$  (cf. \cite{pap-rap}) together with weights at each point of $D$.
The local definition of $\cA_\th$ above leads directly to the definition of a good space of connections on a parahoric bundle.
\vspace{0.2cm}

\noindent{\bf Corollary/conjecture F}.
{\em 
There is a canonical bijection between the set of equivalence classes of triples $(\G, A,E)$  
of connections $A$ on parahoric bundles $E$ for $\G\to\Si$, and {\em filtered $G$-local systems} on $\Si$, i.e. homomorphisms of the fundamental group of $\Si\setminus D$ into $G$ together with a weighted parabolic subgroup of $G$ for each point of $D$ compatible with the local monodromy.
}
}

The layout of this article is as follows.
In section 2 we give basic definitions---this is divided into three parts: reductive groups, loop groups and meromorphic connections.
Section 3 then establishes the main correspondence (Theorem D).
Next section 4 is devoted to quasi-Hamiltonian geometry and establishes Theorem B. Finally section 5 discusses Bruhat--Tits buildings and weighted parahoric subgroups and deduces Corollary E.
Some further directions are mentioned at the end.

\

\noindent{\bf Acknowledgments.}
This research is partially supported by ANR grants 
08-BLAN-0317-01/02 (SEDIGA), 09-JCJC-0102-01 (RepRed).
I would like to thank
O. Biquard, P. Gille,  M. S. Narasimhan, C. Sabbah, C. Simpson and D. Yamakawa.

\section{Basic definitions}\label{sn: defns}

Let $G$ be a connected complex reductive group. 
Let $T\subset G$ be a maximal torus and $B\subset G$ a Borel subgroup containing $T$. Write the Lie algebras as $\lt\subset \lb \subset \g$.
Let $\cR\subset \lt^*$ denote the set of roots and let $\Delta\subset \cR$ denote the simple roots determined by $B$. 
We will identify the roots with characters of $T$ whenever convenient. 
Let $\g_\al\subset \g$ be the root space corresponding to $\al\in \cR$ and  let $U_\al\subset G$ denote the corresponding root group.

Let $X_*(T)$ denote the set of one parameter subgroups $\varphi:\IC^*\to T$ of $T$. Taking the derivative ($\varphi=z^\phi\mapsto \phi$) embeds $X_*(T)$ as a lattice in $\lt$, and we define $\lt_\IR = X_*(T)\otimes_\IZ\IR\subset \lt$, so that $\lt$ is the complexification of the real vector space $\lt_\IR$.

Recall the Jordan decompositions: 1) $X\in \g$ has a unique decomposition $X=X_s+X_n$ with $X_s$ semisimple, $X_n$ nilpotent and $[X_s,X_n]=0$,
2) $g\in G$  has a unique decomposition $g=g_s g_u$ with $g_s$ semisimple, $g_u$ unipotent and $g_sg_u=g_ug_s$.

An element of $X\in \g$ will be said to have {\em real eigenvalues} if its adjoint orbit contains an element whose  semisimple part is in $\lt_\IR$.
Said differently there are a finite number of {\em commuting} one parameter subgroups $\la_i$ such that $X_s = \sum a_i d\la_i$ for real numbers $a_i$.

Recall that the standard parabolic subgroups $P_I\subset G$ are  the subgroups containing $B$.
They are determined by subsets $I$ of the nodes of the Dynkin diagram $\Delta$.
The Lie algebra of $P_I$ is that of $B$ plus the sum of the root spaces $\g_{-\al}$ for positive roots $\al$ which are linear combinations of the elements of $I$.
The parabolic subgroups $P\subset G$ may be characterised as the  subgroups conjugate to a standard parabolic.
The Levi factor of $P$ is the quotient $L=P/U$ of $P$ by the unipotent radical $U=\Rad_u(P)$ of $P$; it is again a connected complex reductive group. One can choose a lifting of $L$ to a subgroup of $P$ (and thus of $G$) and $P$ is isomorphic to the semi-direct product of $L$ and $U$. If $T\subset B\subset P$ then we have a preferred lift $L$ with $T\subset L$, but in general there are many lifts, since we can conjugate the lift $L$ by elements of $P$.

Any semisimple element $\th\in \g$ with real eigenvalues (and in particular any one parameter subgroup) has an associated  parabolic subgroup: 
$$P_\th = \{ g\in G \st z^\th g z^{-\th} \text{ has a limit as $z\to 0$ along any ray}\}\subset G$$
where $z^\th=\exp(\th\log(z))$.
Equivalently $P_{\th}$ is $L\cdot U\subset G$ where the Levi factor $L\subset G$ is the centraliser of $\th$ and $U\subset G$ is the unipotent subgroup whose Lie algebra is the direct sum of the eigenspaces of $\ad_\th \in \End(\g)$ with strictly positive eigenvalues. 
For one-parameter subgroups this notion is used by Mumford \cite{MFK94} p.55.
If we choose $\th$ (or $T$) such that 
$\th\in \lt_\IR$ then $P_\th$ is the group generated by $T$ and the root groups $U_\al$ such that $\al(\th)\ge 0$. (If further $\th$ is in the closed positive Weyl chamber then $P_\th=P_I$ where $I=\{\al\in \Delta \st \al(\th)=0\}$ is the set of walls containing $\th$.) 
Note that $P_{\Ad_h(\th)} = hP_\th h^{-1}$ for any $h\in G$. 

Now let $P\subset G$ be a parabolic subgroup and let $\cC\subset L$ be a conjugacy class in the Levi factor $L$ of $P$.

\begin{lem}\label{lem: conjorbits}
The conjugacy class $\cC\subset L$ uniquely determines a conjugacy class in the Levi factor of any parabolic subgroup of $G$ conjugate to $P$.
\end{lem}
\pf
Given $l\in \cC\subset L$, and $g\in G$ then $glg^{-1}$ projects to an element $h=\pi(glg^{-1})$ of the Levi factor $H$ of the parabolic $Q=g P g^{-1}$
(where $\pi:Q\to H:=Q/\Rad_u(Q)$).
The conjugacy class in $H$ of $h$ is uniquely determined:
Since parabolics are their own normalisers (\cite{Bor91} 11.16)
$Q$ determines $g$ upto left multiplication by an element $q$ of $Q$.
Replacing $g$ by $qg$ only conjugates $h$ by $\pi(q)$. 
Choosing a different $l\in \cC$ corresponds to right multiplication of $g$ by an element $p$ of $P$---this does not change $Q$ so by the  above corresponds to conjugating $h$.
\epf

Similarly an adjoint orbit $\orbit\subset \Lie(L)$ uniquely determines an adjoint orbit of the Lie algebra of the Levi factor of any conjugate parabolic. Similarly also for coadjoint orbits in $\Lie(L)^*$.

Given a parabolic subgroup $P\subset G$, a {\em set of weights} 
for $P$ is an
element  $[\th]$ of the centre of the Lie algebra of the Levi factor $L$ of $P$ such that 1) it is semisimple and has real eigenvalues, and 2) given any lift of $L$ to a subgroup of $P$ the corresponding lift $\th\in \lp\subset \g$ of $[\th]$ determines $P$, i.e. $P_\th = P$.
A {\em weighted parabolic subgroup} is a parabolic subgroup 
$P$ together with a set of weights for $P$. 
More concretely $[\th]$ is a (one-point) adjoint orbit of $L$ 
and so corresponds uniquely to an adjoint orbit of the Levi factor of the standard parabolic $P_I$ conjugate to $P$. Then $[\th]$ just corresponds to a point $\th'$ of the closed Weyl chamber such that $P_I= P_{\th'}$.
Thus if $G$ is semisimple this amounts to choosing a strictly positive real number for each element of $\Delta\setminus I$.

\ppb{
Given a choice of a conjugacy class of parabolic subgroups $P\subset G$ (i.e. the choice of a subset $I\subset \Delta$)
a {\em set of weights for $I$} is an element $\theta\in\lt_\IR$ such that $P_\th = P_I$.
In other words
$\al(\theta)\ge 0$  for all simple roots $\al\in \Delta$ (so $\theta$ is in the closed positive Weyl chamber) and if $\al$ is a simple root then $\al(\theta)>0$ if and only if $\al\in \Delta\setminus I$ (so that $\theta$ determines $I$---indeed $P_I = P_\theta$).
Thus if $G$ is semisimple this amounts to choosing a positive real number for each element of $\Delta\setminus I$.

Thus $\theta$ is in the centre of the Lie algebra of the Levi $L$ of $P_I$, and so is a (one point) adjoint orbit of $L$.
Thus (by the above) $\theta$ uniquely determines an adjoint orbit of the Levi of any  parabolic $P$ conjugate to $P_I$.
(Thus we could give a more intrinsic definition of weights [not using the choice of $T$ and $B$] in terms of certain one point adjoint orbits of Levi factors.)
}

\begin{lem}
A semisimple element $\th\in \g$ with real eigenvalues  determines a set of weights $[\theta]$ for the associated parabolic subgroup $P_\th\subset G$. In general there are many elements $\th$ determining the same pair $(P_\th,[\th])$. 
\end{lem}
\pf Indeed $\th$ determines a Levi decomposition $P_\th=LU$ (with $L$ the centraliser of $\th$) and $\th$ is in the Lie algebra of the centre of $L$, so determines a weight. 
(Less abstractly $\th$ is conjugate to a unique element $\theta'$ of the closed Weyl chamber in $\lt_\IR$.)
Finally it is clear that $\th$ and $g\th g^{-1}$ determine the same pair for any $g\in P_\th$.
\epf 

Let $\IB(G)$ denote the set of weighted parabolic subgroups of $G$.
(This will be discussed in more detail in Section \ref{sn: buildings}.)

\ssn{Background on loop groups}
\nopagebreak

Now we will consider the analogous definitions for the complex (local) loop group.
\ppb{
In $\lt_\IR$ we may consider the affine hyperplane
$$H_{\al,n} = \{ \la\in \lt_\IR \st \al(\la)=n\}$$
for any root $\al$ and integer $n$.
The connected components of $\lt_\IR$ minus all the affine hyperplanes are the alcoves, and there is a unique alcove in the positive Weyl chamber having the origin as a vertex. 
}
We will work with the ring $\cO=\IC\{z\}$ of germs of holomorphic functions (equivalently power series with radius of convergence $>0$)
 and its field of fractions
$\cK=\IC\{\!(z)\!\}=\IC\{z\}[z^{-1}]$.
(The proofs we will give also yield the analogous results for the completions $\wh \cO = \IC\flb z\frb$ and $\wh\cK=\IC(\!(z)\!)$---in fact this case is slightly easier---for simplicity only the completed results were stated in the introduction.)
The convergent local loop group is $LG = G(\cK)$, the group of $\cK$ points of the algebraic group $G$.
The subgroups of $LG$ analogous to parabolic subgroups of $G$ are the {\em parahoric} subgroups of $LG$. 
(Unlike in the finite dimensional case parahoric subgroups are not always self-normalising.)
A basic example of a parahoric subgroup is the subgroup $G(\cO)$ which arises as the group of germs of bundle automorphisms if we choose a local trivialisation of a principal $G$-bundle.
Similarly the Iwahori subgroup 
$$\cI = \{g\in G(\cO)\st g(0)\in B\}$$
and its parabolic generalisations
$$\{g\in G(\cO)\st g(0)\in P\}$$
(where $P\subset G$ is a parabolic subgroup)
arise if we consider parabolic $G$-bundles.
These are also parahoric subgroups of $LG$ but they do not exhaust all the possibilities. 
Indeed, if $G$ is simple, conjugacy classes of parahoric subgroups of $LG$ correspond to proper subsets of the nodes of the  affine Dynkin diagram, 
whereas those above correspond to parabolic subgroups of $G$, i.e. to subsets of the usual Dynkin diagram.
E.g. if $G=E_8$ there are $511$ conjugacy classes of parahoric subgroups of $LG$, of which only $256$ arise from parabolic subgroups of $G$.
On the other hand if $G=\GL_n$ any parahoric subgroup is conjugate to a subgroup arising from a parabolic subgroup of $G$.

The general setup we will need for Theorem D is as follows.
Given an element $\theta\in \lt_\IR$ we will define an associated parahoric subgroup of the loop group.
First $\th$ gives a  grading of the 
Lie algebra $\g$, namely it decomposes as
$$\g = \bigoplus_{\la\in \IR}\g_\la$$
where $\g_\la$ is the $\la$ eigenspace of $\ad_\theta$.
Then for any integer $i$ we may define subspaces
$$\g(i) = \bigoplus_{\la \ge -i} \g_\la\subset \g\qquad\text{and}\qquad
\n(i) = \bigoplus_{\la > -i} \g_\la\subset \g$$
so in particular $\g(0) = \lp_\theta$ is the Lie algebra of the parabolic associated to $\theta$, and $\n(0)$ is its nilradical (and $\g_0$ is its Levi factor).
To emphasise the dependence on $\theta$ we will sometimes write $\g^\theta_\la = \g_\la$ and $\g^\theta(i) = \g(i)$.
Note that the subset
$$\llp_\th:= 
\{X = \sum_{i\in \IZ}X_iz^i \in \g\{\!(z)\!\}) \st X_i\in \g(i) \}$$
is a Lie subalgebra of $L\g=\g\{\!(z)\!\}$.
Said differently $\th$ determines a grading of the vector space $L\g$,
with finite dimensional pieces 
$$L\g(r) = \Bigl\{\sum X_i z^i\in L\g \st X_i\in \g_{\la} \text{ where } \la+i=r\Bigr\}$$ 
for all $r\in \IR$.
Then $\llp_\th$ is the subalgebra of $L\g$ with weights $r\ge 0$.
The weight zero piece will be a  subalgebra
which we will denote as
$$\ll_\th:=L\g(0) = \{ X = \sum X_iz^i \in \llp_\th \st X_i\in \g_{-i} \}.$$

This is finite dimensional  and  in fact reductive.
We view $\ll_\th$  as the Levi factor of $\llp_\th$. 
By setting $z=1$ there is an embedding
$$\iota:\ll_\th\hookrightarrow\g.$$
Let $\wh H_\th = C_G(e^{2\pi i \th})$ be the centraliser in $G$ of 
$e^{2\pi i \th}$, and let $\lh_\th\subset \g$ be its Lie algebra.
Then the image $\iota(\ll_\th)$ is $\lh_\th$.
(Note that $\lh_\th$ is not necessarily isomorphic to a Levi factor of a parabolic subalgebra of $\g$---for example for simple $\g$, $\lh_\th$ could be the Lie algebra determined by any proper subset of the nodes of the affine Dynkin diagram of $\g$, so may be a proper semisimple subalgebra of the same rank, such as $\sl_3\subset \g_2$, as in Borel--De Siebenthal theory.)
More generally we may consider the subgroup
$$\wh \cL_\th = \{ z^{-\th} h z^\th \st h\in \wh H_\th\}\subset LG$$
of  $LG$ (this is indeed well defined since $h$ commutes with the monodromy of $z^\th = \exp(\th \log(z))$). By setting $z=1$ we see 
$\wh \cL_\th$ is isomorphic to $\wh H_\th$, and $\iota$ is the corresponding map on the level of Lie algebras.
Let $H_\th$ denote the identity component of $\wh H_\th$ and let $\cL_\th\subset \wh \cL_\th$ denote the corresponding subgroup of the loop group.
Thus the Lie algebra of $\wh \cL_\th$ and $\cL_\th$ is $\ll_\th$.

The extended parahoric subgroup determined by $\th$ is the subgroup 
$$\wh \cP_\th = \{ g\in LG \st z^\th g z^{-\th} \text{ has a limit as $z\to 0$ along any ray}\}.$$
This definition is perhaps best understood by thinking in terms of a faithful representation, whence $\th$ is a diagonal matrix and we can see explicitly what the condition means in terms of matrix entries.
Alternatively one can work with the Bruhat decomposition, and show that $\wh \cP_\th$ is generated by 1) elements of $\wh \cL_\th$, 2) elements of the form $\exp(X z^i)$ with $X\in \g_\al$ such that $\al(\th)+i > 0$
(or $X\in \lt$ and $i>0$)
and 3)  elements of the form $\exp(Y(z))$ with $Y\in z^N\g\{ z \}$
with $N$ a sufficiently large integer (so that $Y\in \llp_\th$).   
Heuristically the Lie algebra of $\wh \cP_\th$ is $\llp_\th$.
This has Levi subgroup
$\wh \cL_\th$ and pro-unipotent radical 
$$ \cU_\th = \{ g\in LG \st z^\th g z^{-\th} \text{ tends to $1$ as $z\to 0$ along any ray}\}.$$
(which has Lie algebra the part of $L\g$ of weight $>0$, and is generated by elements just of type 2) and 3) above).
The group $\wh \cP_\th$ is the semidirect product of $\wh \cL_\th$ and $\cU_\th$.

The parahoric subgroup associated to $\th$ is the group
generated by  $\cU_\th$ and the connected group $\cL_\th$:
$$\cP_\th = \cL_\th \cdot \cU_\th \subset \wh \cP_\th.$$
This is a normal subgroup of $\wh \cP_\th$ and the quotient 
$\wh \cP_\th/\cP_\th \cong \wh H_\th/H_\th$ is finite.

\ppb{
The ``nilpotent radical'' of $\llp_\th$ is the subalgebra of elements of positive weight
$$\lln_\th = \{ X = \sum X_iz^i \in \llp_\th \st X_i\in \n(i) \}.$$
Let $\cP_\theta\subset LG$ be the group corresponding to $\llp_\th$.
It is the semidirect product of $H_\theta$ and the group $\cU_\th$ corresponding to $\lln_\th$ (which we think of as the unipotent radical of $\cP_\th$).
Given an integer $N$ such that $\g(N) = \g$ (so that 
$z^N\g(\cO)\subset \llp_\th$)  we may consider the group $\G(N)$ of germs at $0$ of holomorphic maps from neighbourhoods of 
$0\in \IC$ to $G$ which are tangent to order $N$ to $1\in G$ at $0$. Any element of $\cP_\th$ may then be written as a finite product of terms of  $H_\th$, of $\G(N)$ and of the form $\exp(z^iX)$ (with 
$X\in \g_\la$ for some $\la$ with $\la+i>0$).

Let $N(\cP_\th)$ denote the normaliser of $\cP_\th$ in $LG$.

Question: (for $G$ ss sc) is this the set of $g\in LG$ stabilising $\llp_\th \subset Lg$? ???

Facts: If $G$ is simply connected then $N(\cP_\th)=\cP_\th$. In general 
$N(\cP_\th)$ is obtained by replacing $H_\th$ by a disconnected group with the same identity component.

??? the centraliser in LG? of its center
??? isomorphic to the centraliser in $G$ of $\exp(2\pi i \th)$?
}

\ssn{Germs of meromorphic connections}
\nopagebreak

Choose $\th\in \lt_\IR$ and let $\cP_\th$ be the corresponding parahoric subgroup with Lie algebra $\llp_\th$.
Then we may consider the space $\cA=\g(\cK) dz$ of meromorphic connections (on the trivial 
$G$-bundle over the disc) and the subspace
$$\cA_\th = \llp_\th \frac{dz}{z}.$$

Thus if $\th=0$ this is just the space of logarithmic connections.
If $\th$ is small, these are the logarithmic connections with residue in the Lie algebra $\lp_\th$ of $P_\th$, as occurs in the case of parabolic bundles.
(Parabolic $G$-bundles are studied for example in  \cite{telwood-par}, in the case where $G$ is simple and simply-connected, and the weights are small and {rational}.)
In general elements of $\cA_\th$ will have poles of order greater than one, but we will see in the course of the proof of Theorem \ref{thm: main corresp} below that they always have regular singularities: fundamental solutions have at most polynomial growth at zero.  
They should perhaps be viewed as the right notion of ``logarithmic parahoric connections'' (as the pole is of order one greater than that permitted by the parahoric structure) but this term is cumbersome and possibly confusing.
We will call them tame parahoric connections (although perhaps {\em ``logahoric''} is simplest).

\begin{lem} \label{lem: wP preserves Ath}
The natural (gauge) action of $\wh \cP_\th$ on $\cA$ preserves $\cA_\th$.
\end{lem}
\pf
Given $g\in LG$ and a connection $A\in \cA$, 
the gauge action of $g$ on $A$ is
$g[A] := \Ad_g(A) + (dg)g^{-1}$
where for any $g\in LG$ 
we define the  $\g$-valued  meromorphic one-form 
(on a neighbourhood of $0\in \IC$) 
$$(dg)g^{-1}:=g^*(\bar\Theta) $$
where $\bar\Theta\in \Omega^1(G,\g)$ is the right-invariant Maurer--Cartan form on $G$. (The sign conventions used here are as in \cite{bafi}.)
Thus it is sufficient to check that $(dg)g^{-1}\in \cA_\th$ for any $g\in \wh \cP_\th$. First if $g=\exp(X z^i)$ for $X\in \g_\al$ with $\al(\th)+i \ge 0$ then $(dg)g^{-1} = (i X z^i)dz/z\in \cA_\th$. Second if $g=\exp(X(z))$ with $X\in z^N\g\{ z \}$ for $N$ a sufficiently large integer
again we have $(dg)g^{-1}\in \cA_\th$. Such elements generate $\cP_\th$ so it follows that  $\cP_\th$ preserves $\cA_\th$ (since 
$d(gh)(gh)^{-1} = \Ad_g(dhh^{-1}) + dgg^{-1}$).
Finally we must check $\wh \cL_\th$ preserves $\cA_\th$ (since $\wh \cP_\th$ is generated by this and $\cP_\th$).
But if $g=z^{-\th} h z^{\th}\in \wh \cL_\th$ then one finds $dgg^{-1} = (\Ad_g(\th) -\th)dz/z$ and this will be in $\cA_\th$ if 
$\Ad_{z^\th} (\Ad_g(\th) -\th)$ has a limit as $z\to 0$ along any ray.
But it does have a limit, since it is constant and equals $\Ad_h(\th) -\th$.
\epf

Note that the gauge action may also be interpreted as the (level one) coadjoint action of a central extension of $LG$, although we will not need this interpretation here.
A closer examination of the action of $\wh \cL_\th$ on the  weight zero piece
$\cA_\th(0)=L\g(0) dz/z$ of $\cA_\th$ yields the following.

\begin{lem}\label{lem: local orbits}
The map $A\mapsto z^\th[A]$ is well defined on the weight zero piece 
$\cA_\th(0)$ of $\cA_\th$ and provides an isomorphism 
$$\cA_\th(0)\cong \lh_\th\frac{dz}{z}\subset \g \frac{dz}{z}$$
which is equivariant with respect to the gauge action of $\wh \cL_\th$ on 
$\cA_\th(0)$ and the adjoint action of $\wh H_\th$ on $\lh_\th$.
\end{lem}
\pf $\cA_\th(0)$ is just the set of elements $A=\sum A_iz^i dz/z$ with 
$A_i\in \g$ in the $-i$ eigenspace of $\ad_\th$.
Thus $z^\th[A]=(\th+\sum A_i) {dz}/{z}\in \g {dz}/{z}$, and $B:=\th+\sum A_i$ is just an arbitrary element of $\lh_\th$
(i.e. an element with components only in the integer eigenspaces of $\ad_\th$). Clearly if $h\in \wh H_\th$
and  $g=z^{-\th} h z^{\th}\in \wh \cL_\th\subset LG$
is the corresponding element of the loop group
then
$$  z^\th \left[g[A]\right ] = h [B dz/z] =  \Ad_h(B) \frac{dz}{z}$$
so we have the desired equivariance.
\epf

\ppb{

Consider pairs $(A,\cP)$ where $\cP\subset LG$ is a parahoric subgroup and $A$ is a compatible connection, namely an element of
$$ \cA = \llp \frac{dz}{z},$$
(where $\llp$ is the Lie algebra of $\cP$) on which $\cP$ and its normaliser $\wh \cP$ act by gauge transformations.
Then $LG$ acts on the set of pairs $(A,\cP)$ and we may consider the set of orbits. This amounts to computing the quotient $\cA/\wh\cP$ for all
conjugacy classes of parahoric subgroups. Thus it is sufficient to fix $\th\in \lt_\IR$ and study  $\cA_\th/\wh\cP_\th$ where $\cA_\th = \llp_\th \frac{dz}{z}.$

Given $g\in LG$ and a connection $A\in L\g dz/z$, the gauge action of $g$ on $A$ is
$$g[A] := \Ad_g(A) + (dg)g^{-1}.$$
where for any $g\in LG$ 
we define the  $\g$-valued  meromorphic one-form 
(on a neighbourhood of $0\in \IC$) 
$$(dg)g^{-1}:=g^*(\bar\Theta) $$
where $\bar\Theta\in \Omega^1(G,\g)$ is the right-invariant Maurer--Cartan form on $G$.
This may also be interpreted as the (level one) coadjoint action of a central extension of $LG$, although we will not need this interpretation here.

To start with let us examine the weight zero piece of $\cA_\th$.

\begin{lem}\label{lem: local orbits}
The map $A\mapsto z^\th[A]$ is well defined on the weight zero piece 
$\cA_\th(0)=\lh_\th dz/z$ of $\cA_\th$ and provides an isomorphism 
$$\cA_\th(0)\cong \iota(\lh_\th)\frac{dz}{z}\subset \g \frac{dz}{z}$$
which is equivariant with respect to the gauge action of $\wh H_\th$ on the left hand side and the adjoint action of $\iota(\wh H_\th)=C_G(e^{2\pi i \th})$ on the right.
\end{lem}
\pf $\cA_\th(0)$ is just the set of elements $A=\sum A_iz^i dz/z$ with 
$A_i\in \g$ in the $-i$ eigenspace of $\ad_\th$.
Thus $z^\th[A]=\sum A_i {dz}/{z}\in \g {dz}/{z}$, and $B:=\sum A_i$ is just an arbitrary element of the Lie algebra of the centraliser of 
$e^{2 \pi i \th}$ (i.e. an element with components only in the integer eigenspaces of $\ad_\th$). Clearly if $h\in C_G(e^{2\pi i \th})$ then
then $h=\iota(g)$ where $g=z^{-\th} h z^{\th}\in \wh H_\th\subset LG$
and 
$$  z^\th \left[g[A]\right ] = h [B dz/z] =  \Ad_h(B) \frac{dz}{z}$$
so we have the desired equivariance.

\epf

Suppose we have a parahoric subgroup $\cP_\th\subset LG$
In a local trivialisation, a meromorphic connection  compatible with the parahoric structure is an element of 
$$\cA_\th := \llp_\th\frac{dz}{z}$$
and two connections $A,B$ are isomorphic if they are in the same orbit of $\cP_\th$ acting by gauge transformations., i.e. $A = g[B]$ for some $g\in \cP_\th$ where 
$$g[B] := \Ad_g(B) + (dg)g^{-1}.$$
For any $g\in LG$ we define 
the $\g$-valued meromorphic one-form 
(on a neighbourhood of $0\in \IC$) 
$$(dg)g^{-1}:=g^*(\bar\Theta) $$
where $\bar\Theta\in \Omega^1(G,\g)$ is the right-invariant Maurer--Cartan form on $G$.
Note that if $g\in \cP_\th$ then $(dg)g^{-1}\in \llp_\th\frac{dz}{z}$.
(This follows since it is true for elements of $\G(N)$ for large $N$, and for elements of the form $\exp(Xz^i)$ for $X\in \g_\al$ with $\al(\th)+i\ge 0$, and these generate $\cP_\th$.)

???in general look at action of normaliser of $\cP_\th$.

Note that if $\al(\theta)\le 1$ for all roots $\al$ then $\g(i)=\g$ for all $i\ge 1$, and so $\llp_\theta$ contains the subset of $\g\{z\}$ with constant term in $\lp_\theta$.
Thus in this case $\cP_\th$ contains the standard Iwahori subgroup (provided we choose a Borel subalgebra in $\lp_\th$).
If further $\al(\theta) < 1$ for all roots $\al$ then  $\g(i)=\{0\}$ for $i<0$ so $\llp_\th\subset \g\{z\}$ and $\cA_\th$ consists of logarithmic connections with residue fixed to preserve a parabolic structure (this is the ``classical case'' of connection on parabolic $G$-bundles.
If $\th=0$ then we reduce to the case of logarithmic connections.

??? repitition here
}
\section{Main correspondence} 

Having now covered the background definitions we can move on to the main result.
Fix elements $\th, \tau \in \lt_\IR$ and $\si\in \sqrt{-1}\lt_\IR$ and set $\phi=\th+\tau$.
Choose a nilpotent element $n\in \g$ commuting with $\phi$ and $\si$ and such that 
$\Ad_t(n) = n$ where $t:= \exp(2\pi i \tau)\in G$.
(This means there is a (finite) decomposition $n=\sum a_i$ with $[\tau,a_i]=ia_i$ for $i\in \IZ$.)

As above $\th$ determines a
space $\cA_\th$  of $\th$ parahoric connections and an 
extended parahoric subgroup 
$\wh \cP_\th\subset LG$ with Levi subgroup $\wh \cL_\th$.
Moreover $\th$ determines an isomorphism $\wh \cL_\th\cong \wh H_\th:= C_G(e^{2\pi i \th})\subset G$. The corresponding Lie algebras are denoted $\ll_\th\cong \lh_\th\subset \g$.

Let $\orbit\subset \lh_\th$ be the adjoint orbit (under the possibly disconnected group $\wh H_\th$) of 
the element
$$\phi + \si + n \in \lh_\th.$$ 
This corresponds to the $\wh \cL_\th$ orbit in $\cA_\th(0)$ containing the element 
\beq\label{eq: DR local data}
(\tau + \si + \sum a_i z^i)\frac{dz}{z}.\eeq 

Also $\phi$ determines a parabolic subgroup $P_\phi\subset G$ and a weight $[\phi]$ for $P_\phi$. 
Let $L$ be the centraliser of $\phi$ (a Levi subgroup of $P_\phi$).
By construction $\tau,\si$ and $n$ commute with $\phi$, so are in the Lie algebra of $L$. Then we define $\cC\subset L$ to be the conjugacy class
containing the element
$$\exp(2\pi i (\tau+\si))\exp(2\pi i n )\in L.$$

Note that $\cC$ is not necessarily an exponential conjugacy class, since $n$ and $\tau$ in general do not commute. (This was one of our motivations for considering more general objects than  logarithmic or parabolic connections.) 

\begin{lem} \label{lem: orbit idpce}
The triple $([\phi], P_\phi, \cC)$ is uniquely determined upto conjugacy by $\th$ and the orbit $\orbit\subset \lh_\th$. Moreover any such triple $([\phi], P_\phi, \cC)$ arises in this way (upto conjugacy).
\end{lem}
\pf
Any other element of $\orbit$ will be of the form $\Ad_g(\phi+\si+n)$ with $g\in \wh H_\th$. Suppose we choose $g$ so that the semisimple part $\Ad_g(\phi+\si)$ is in $\lt$.
This yields new choices $\phi' =\Ad_g\phi , \tau' = \phi' - \th, \si' = \Ad_g\si , n' = \Ad_g n$ and we should check that
$([\phi'], P_{\phi'}, \cC')$ is conjugate to $([\phi], P_\phi,\cC)$.
But this follows from the fact that
$\exp(2\pi i \tau') = g\exp(2\pi i \tau )g^{-1}$, as $g$ commutes with $\exp(2\pi i \th)$.
The fact that all such triples arise follows immediately from the multiplicative 
Jordan decomposition.
\epf

Note that, given $\phi=\th+\tau$ and $\si$, the precise correspondence between the adjoint orbits $\orbit\subset \lh_\th$ and the conjugacy classes $\cC\subset L=C_G(\phi)$ rests on the identification
\begin{align*}
\{X\in \lh_\th \st [X,\phi]=[X,\si]=0 \} 
&= \{X\in \g \st [X,\phi]=[X,\si]=0, \Ad_{e^{2\pi i \th}}X=X \} \\
&= \{X\in \g \st [X,\phi]=[X,\si]=0, \Ad_{e^{2\pi i \tau}}X=X \} \\
&=  \{X\in \Lie(L) \st [X,\si]=0, \Ad_{e^{2\pi i \tau}}X=X \}
\end{align*}
since $n$ is a nilpotent element of this (reductive) Lie algebra.
\ppb{

Since $\lt$ is a Cartan subalgebra of $\lh_\th$ (which is reductive) we may choose an element of the form 
$$\phi+\si + n$$
in $\orbit$ with $\phi\in \lt_\IR, \si\in \sqrt{-1}\lt_\IR$ and 
$n$ a nilpotent element commuting with $\phi$ and $\si$. 
In turn we may write $n=\sum a_i$ with $a_i$ in the $-i$ eigenspace of $\ad_\th$.
Then we may set $\tau = \phi-\th$ and consider the element
$$X\frac{dz}{z}\in \cA_\th\qquad\text{where}\qquad X = \tau+\si + \sum a_i z^i.$$
By construction this is in the weight zero piece $\lh_\th dz/z$ 
of $\cA_\th$. 
We also have $[\tau, a_i]=ia_i$ for all $i$.

Note that the group $H_\th$ acts on $\lh_\th dz/z$ by gauge transformations, and the orbits of this action match up precisely (in the above fashion) with the adjoint orbits of $\iota(H_\th)\subset G$ on $\iota(\lh_\th)\subset \g$.
}

We will say that a connection $A\in \cA_\th$ ``lies over $\orbit$'' if 
its weight zero component is in the $\wh\cL_\th$ orbit 
corresponding to $\orbit$.
Similarly if $P\subset G$ is a parabolic subgroup conjugate to $P_\phi$ we will say $M\in P$ ``lies over $\cC$'' if $\pi(M)\in \cC$, where $\pi$ is the canonical projection from $P$ onto its Levi factor (and we transfer $\cC$ from $L\subset P_\phi$ as in Lemma \ref{lem: conjorbits}).

The main statement (Theorem D of the introduction) is then:

\begin{thm}\label{thm: main corresp}
There is a canonical bijection between the $\wh \cP_\th$ orbits of 
tame parahoric connections in $\cA_\th$ lying over $\orbit$ and conjugacy classes of pairs $(M,P)$ with $P\subset G$ a parabolic subgroup conjugate to $P_\phi$ and $M\in P$ an element lying over $\cC$. 
\end{thm}

\ppb{
Clearly setting $\th=0$ (so $\phi=\tau$) yields the result of the introduction.

???really should phrase in terms of $LG$ orbits of pairs $(A,\cP_\th)$
such that $A\in \llp_\th dz/z$---this means taking into account the normaliser of $\cP_\th$ in the non sc case.)
}
\pf
To start we will explain how to put such connections in a simpler form.
Suppose $A\in \cA_\th$ lies over $\orbit$.
First we may do a gauge transform by an element of $\wh \cL_\th$ 
so the weight zero component of $A$ equals $A(0):= (\tau+\si+\sum a_iz^i) dz/z$.
Then we claim we can do a gauge transformation by an element $g$ of $\cU_\th$ such that $g[A]$ is normalised in the following way:
\beq
g[A] = (\tau + \si + \sum_{i\in \IZ} A_i z^i)\frac{dz}{z}
\eeq
with  each $A_i\in \g(i)$ and
$$[\tau,A_i]=iA_i,\qquad [\si,A_i] = 0$$
for all $i\in \IZ$ (and $a_i$ is the component of $A_i$ in the $-i$ eigenspace of $\ad_\th$). This implies only finitely many of the $A_i$ are nonzero.

To prove the claim we extend the usual argument in the logarithmic case   (cf. \cite{BV83}) as follows.
Let $0=r_0<r_1<\cdots$ be the sequence of positive real numbers such that 
$L\g(r_i)\ne 0$. Suppose inductively that the piece of $A$ in $L\g(r_i)\frac{dz}{z}$ has been normalised for $0\le i < k$.
Then we claim we may choose $X(k)\in L\g(r_k)$ so that the piece of $g_k[A]$ in  $L\g(r_i)\frac{dz}{z}$ is normalised for $0\le i \le k$,
where $g_k = \exp(X(k))\in \cU_\th$.
Indeed $g_k[A]$ will equal $A$ up to $L\g(r_{k-1})\frac{dz}{z}$ and will have subsequent coefficient
\beq\label{eq:next coeff}
A(k) + [X(k),A(0)] + z\frac{d}{dz}X(k)\in L\g(r_{k})
\eeq
where $A = \sum_{i\ge 0}A(i)\frac{dz}{z}$ with $A(i)\in L\g(r_i)$.
Thus ideally we would like to choose $X(k)$ such that this was zero, i.e.
$\left(\ad_{A(0)} - z\frac{d}{dz}\right) X(k) = A(k)$.
This is not always possible, but we can make the difference small, as follows.
Note that $\ad_{A(0)}$ restricts to a linear operator on the finite dimensional vector space $L\g(r_k)$ and so we may decompose $L\g(r_k)$
into its generalised eigenspaces. 
Since $\tau+\si$ is the semisimple part of $A(0)$, these  generalised eigenspaces are just the eigenspaces of the semisimple operator 
$\ad_{\tau+\si}$. On the other hand we also have 
$z\frac{d}{dz}\in \End(L\g(r_k))$ preserving this eigenspace decomposition
(as it commutes with $\ad_{\tau+\si}$)
and having only integral eigenvalues (mapping $xz^i$ to $ixz^i$).
Thus if $a_{i\mu}z^i$ (resp. $x_{i\mu}z^i$) 
is the component of  $A(k)$ (resp. $X(k)$) in the $\mu$-eigenspace of 
$\ad_{\tau+\si}$ (and the $i$-eigenspace of $z\frac{d}{dz}$) then we may define $X(k)$ by setting $x_{i\mu}=0$ if $\mu=i$ and
$$ x_{i\mu}z^i = \left(\ad_{A(0)} - z\frac{d}{dz}\right)^{-1}( a_{i\mu}z^i)$$
if $i\ne \mu$, since then the operator $\ad_{A(0)} - z\frac{d}{dz}$ will be invertible on the corresponding joint eigenspace.
If $X(k)$ is defined in this way we thus find that the next coefficient \eqref{eq:next coeff} is the sum of the components of $A(k)$ with $i=\mu$, i.e. $\sum_i a_{i i}z^i\in L\g(r_k)$ (noting that $\th$ commutes with $\tau+\si$ and $d/dz$).
But this just means it is normalised: $[\tau,a_{ii}]=ia_{ii}, [\si,a_{ii}]=0$.
Thus inductively we may construct a formal transformation $g=\cdots g_3g_2g_1$ in the completion of $\cU_\th$ converting $A$ into normal form. 
To conclude that this transformation is actually convergent we need to check that 

\begin{lem}
Any connection $A\in \cA_\th$ is regular singular.
\end{lem}
\pf
Choose a faithful representation of $G$ and work in this representation.
Thus $\tau$ is a real diagonal matrix and we may choose 
a diagonal matrix $\la$ with integral eigenvalues such that the diagonal entries of $\tau-\la$ are in $[0,1)$. Let $\varphi = z^\la$ be the corresponding one parameter subgroup.
We may then choose $k$ sufficiently large such that the 
convergent meromorphic gauge transformation $\varphi^{-1} g_k\cdots g_2g_1$ converts $A$ into a logarithmic connection.
\epf

Thus both the original connection and the resulting connection are convergent connections with regular singularities, and it follows that $g$ is actually convergent and in $\cU_\th$ (in effect we constructed an $\wh \cO$ point of a group scheme and then deduced it is actually an $\cO$ point, i.e. in $\cU_\th$).

Having completed the normalisation now  set $N=\sum A_i$, $R = \si + N$, 
$M_s = \exp(2\pi i (\tau+\si)), M_u=\exp(2\pi i N), M= M_sM_u\in G$.

Then $M_sM_u$ is the Jordan decomposition of $M$
and the connection $A$ has monodromy $M$.
Indeed $A$ (after normalisation) equals $z^\tau[Rdz/z]$, which has fundamental solution $z^\tau z^R$.
By construction $M\in P_\phi$, so we have attached a pair $(M,P)$ to the original data (with $P=P_\phi$).
If $L$ is the  Levi factor of $P$ then the image of $M$ in $L$ is 
$$\pi(M) = \exp(2\pi \sqrt{-1} (\tau+\si))\exp\left(2\pi \sqrt{-1} \sum a_i\right)$$
where $a_i$ is the component of $A_i$ in the $-i$ eigenspace of $\ad_\theta$ (the component commuting with $\phi$). 
It follows that  $\pi(M) \in \cC$.

{\bf Surjectivity.\ } 
To give the inverse construction we proceed as follows.
Suppose we have $(M,P)$ with $M\in P$, $P$ of type $\phi$ and $\pi(M)\in \cC\subset L(P)$.
Then we may conjugate by $G$ so that $P=P_\phi$ and $M_s=\exp(2\pi i (\tau+\si))\in T$, since $\pi(M)\in \cC$. (Here $M_s$ is the semisimple part of $M$.)
Then we may write $M=\exp(2\pi i \tau)\exp(2\pi i R)$ for a unique element $R=\si+N\in \g$ with $N$ nilpotent (and $\tau, \si$ as fixed above).
Moreover $N$ commutes with $\si$ and $\Ad_{t}N=N$ where $t:=\exp(2\pi i \tau)\in T$, but $N$ does not necessarily commute with $\tau$ itself. This implies that there is a unique decomposition 
$N = \sum_{i\in \IZ} A_i$ with $A_i\in \g$ such that $[\tau,A_i]=iA_i$.
(If $N$ has components in any other eigenspace of $\ad_\tau$ then one will not have $\Ad_t(N)=N$.)
On the other hand, since $M\in P_\phi$ we have $N\in \lp_\phi$ 
and so $N$ only has components in the positive weight spaces of $\phi$.
Now since $\th=\phi-\tau$
the connection $A:=z^\tau[Rdz/z] = (\tau+\si+\sum A_i z^i)dz/z$ is in $\cA_\th$. 

The component of $A$ in the weight zero component 
$\cA_\th(0)$ is  $(\tau+\si+\sum a_iz^i)dz/z$ 
where $a_i$ is the component of $A_i$ commuting with $\phi$ (i.e. weight $-i$ for $\th$).
This is determined by $\cC$ and lies over $\orbit$ by construction.
Thus $(M,P)$ is the data attached to the connection $A$.

The main lemma we need for the rest of the proof is the following.
\begin{lem}\label{lem: main lem}
Suppose $C\in G$. Then $z^\tau C z^{-\tau}$ is in $\wh \cP_\th$ if and only if a) $C\in C_G(t)$ where $t=e^{2\pi i\tau}$, and b) $C\in P_\phi$.
\end{lem}
\pf
Condition a) holds if and only if $p:=z^\tau C z^{-\tau}$ is in $LG$ (since that is the condition for it to have no monodromy).
Then by definition  $p\in \wh\cP_\th$ if and only if 
$z^\th p z^{-\th}$ has a limit as $z\to 0$ along any ray, i.e. if 
$z^\phi C z^{-\phi}$ has a limit as $z\to 0$. But this is just the condition for $C\in P_\phi$. 
\epf

{\bf Well-defined on orbits.\ }
Next we will check that if two connections are in the same orbit then their data $(M,P)$ are conjugate.
Recall we have fixed $\th, \tau\in \lt_\IR$ and set $\phi = \tau+\th$.
Suppose $A, B\in \cA_\th$ are related by $g\in\wh\cP_\th$.
Without loss of generality we may assume $A,B$ are both normalised.
Thus they have fundamental solutions 
$\Phi_A= z^\tau z^R$ and 
$\Phi_B= z^\tau z^{R_1}$. 
Their monodromies are
$M(A)=te^{2\pi i R}$ and 
$M(B) = te^{2\pi i R_1}$ (both in $P_\phi$) where $t:= e^{2\pi i \tau}$.
The hypothesis means that $\Phi_A = g \Phi_B C$ for some $C\in G$.
This implies $M(A) = C^{-1}M(B)C$, so the monodromies are conjugate, but we must show that $C\in P_\phi$.
It follows (from $M(A) = C^{-1}M(B)C$) that $C$ commutes with $t$ and that $R_1 = \Ad_C(R)$. Thus the identity $\Phi_A = g \Phi_B C$ simplifies to $z^\tau = g(z)z^\tau C$, so that $z^\tau C z^{-\tau} = g^{-1} \in \wh \cP_\th$.
Thus by Lemma \ref{lem: main lem}, $C\in P_\phi$ as desired.

{\bf Injectivity.\ }
Now suppose $A, B\in \cA_\th$ both lie over $\orbit$ and yield data conjugate to $(M,P)$. We will show they are gauge equivalent by an element of $\wh \cP_\th$.
Without loss of generality we may assume $P=P_\phi$.
Thus they have fundamental solutions 
$\Phi_A= f(z)z^\tau z^R$ and 
$\Phi_B=h(z)z^\tau z^{R_1}$ respectively for some 
$f,h\in \wh\cP_\th$, and they have monodromy 
$M(A)=te^{2\pi i R}$ and 
$M(B) = te^{2\pi i R_1}$ where $t:= e^{2\pi i \tau}$.
By assumption $M(B) = CM(A)C^{-1}$ for some $C\in P_\phi=N_G(P_\phi)$.
Thus $M(B) = e^{2\pi i \tau_1}e^{2\pi i \Ad_C(R)}$ where 
$\tau_1 = \Ad_C(\tau)$. 
Comparing the semisimple parts of the two expressions for $M(B)$ we deduce $C$ commutes with $t$ (so that $t=e^{2\pi i \tau_1}$) and also that 
$\Ad_C(\si)=\si$.
Using e.g. the Iwasawa decomposition it follows that $R_1 = \Ad_C(R)$.
Thus $B$ also has fundamental solution 
$h(z)z^{\tau}C z^R=\Phi_BC$.
Thus it is sufficient to prove that $$p:=z^{\tau}Cz^{-\tau} \text{ is in $\wh \cP_\th$}$$
 (since then $h(z)z^{\tau}C z^R \Phi_A^{-1} = hz^{\tau}Cz^{-\tau}f^{-1}$ will be an element of $\wh \cP_\th$ relating $A$ and $B$). 
But by Lemma \ref{lem: main lem} this is now immediate.
\epf

This establishes the main correspondence. For small weights $\th$ this reduces to the parabolic statement in Theorem C (the connected centraliser condition ensures $\wh \cP_\th = \cP_\th$; it is the group $\cP_\th$ that appears in the local moduli of parabolic bundles).
For $\th=0$ one obtains the logarithmic statement (Theorem A).

\begin{rmk}
Note it follows from the proof that the stabiliser in $\wh \cP_\th$ of a connection in $\cA_\th$ is isomorphic to the centraliser in $P_\phi$ of the monodromy $M$. Indeed for a connection in 
normal form this correspondence is given by 
$C\in C_{P_\phi}(M) \leftrightarrow z^\tau C z^{-\tau}$,
and in general one conjugates by any transformation putting the connection in normal form.  
\end{rmk}

\begin{rmk}
Analogously to \cite{Sim-hboncc} one may define the notion of ``filtered $G$-local system'' on a smooth punctured Riemann surface $U$, to be a
$G$-local system $\IL$ on $U$ together with (on a small punctured disc $\Delta_i$ around the $i$th puncture, for each $i$)  a $P$-local system $\IL_i$ (for some  weighted parabolic $P\subset G$) such that  the restriction of $\IL$ to $\Delta_i$ is the $G$-local system 
$\IL_i\times_P G$ associated to $\IL_i$.
If $U$ is a punctured disc, and we choose a basepoint in $U$, then specifying a filtered $G$-local system is the same as specifying the data $(M,P,[\phi])$ in our correspondence (so the correspondence could be restated more intrinsically in terms of filtered $G$-local systems).
\end{rmk}

\begin{rmk}
Another motivation for studying such ``enriched'' (or ``exact'') Riemann--Hilbert correspondences is related to {\em isomonodromic deformations}.
For example such {``monodromy preserving'' deformations} of a nonresonant logarithmic connection $A=\sum A_i dz/(z-a_i)$ on the trivial bundle on $\IP^1$ 
are governed by Schlesinger's equations: 
$dA_i = -\sum_{j\ne i} [A_i,A_j] d\log(a_i-a_j)$.
These are the deformations which preserve the conjugacy class of the monodromy representation of $A$.
Of course these equations make sense for any residues, and one may ask what exactly is preserved by Schlesinger's equations in the resonant case?\footnote{
Recall a logarithmic connection on a $G$-bundle $E$ is {\em resonant} if the residue of the induced connection on the associated vector bundle $\Ad(E)$  has an eigenvalue in $\IZ\setminus\{0\}$.} 
The answer (which is clear from \cite{Mal-imd1long}, or may be extracted from \cite{bolib-fimds}) is that the monodromy representation and the filtrations are preserved (upto overall conjugacy).
This now extends immediately to arbitrary $G$.
Such ``resonant'' deformations are important since for example soliton solutions arise as such (when one has a further irregular singularity at infinity). 
\end{rmk}

\ppb{
\begin{rmk}
Another motivation for studying such ``enriched'' (or ``exact'') Riemann--Hilbert correspondences was to understand the data that is preserved under {\em isomonodromic deformations} of meromorphic connections on Riemann surfaces.  (For this it is convenient to define a meromorphic connection to include a $G$-bundle on a compact Riemann surface.)
In the De\! Rham viewpoint one 
deforms the connection so it fits into an integrable absolute connection.
In the nonresonant logarithmic case this amounts to preserving the monodromy representation, and on $\IP^1$ such deformations are governed by Schlesinger's equations.
In the (generic) nonresonant irregular case it was recognised by Jimbo--Miwa--Ueno \cite{JMU81} that these integrable deformations correspond to preserving the monodromy and the Stokes data (the Stokes data is the natural generalisation of the monodromy representation and so the term ``isomonodromic'' is kept).
[The exact definition of ``preserved'' is best understood in terms of the monodromy data forming a nonlinear local system over the parameter space---the Betti viewpoint in \cite{smid} section 7.] 
This begs the question as to what is preserved in such an integrable deformation of a resonant connection (for example evolving a Fuchsian system according to Schlesinger's equations): The answer is that the filtrations are preserved (as well as the monodromy/Stokes data). (Such deformations are in fact important---for example soliton solutions are such a deformation \cite{JMU81} 6.1).\footnote{ 
Again we vouch for applying the not technically correct term {\em isomonodromic} to such deformations, although {\em isofistomic} is another possiblilty.}
\end{rmk}
}

\section{Quasi-Hamiltonian spaces}

Fix a connected complex reductive group $G$ and a parabolic subgroup $P_0\subset G$. %
Let $\cC\subset L$ be a conjugacy class of the Levi factor $L$ of $P_0$.
Let $\IP\cong G/P_0$ be the variety of parabolic subgroups of $G$ conjugate to $P_0$.

The aim of this section is to prove Theorem B, that 
the set $\wh \cC$ of pairs $(g,P)\in G\times \IP$ with $g\in P$ and $\pi(g)\in \cC$, is a quasi-Hamiltonian $G$ space with $G$-valued moment map given by 
$(g, P)\mapsto g$.

Recall (cf. \cite{AMM}) that a complex manifold $M$ is a 
{\em complex quasi-Hamiltonian $G$-space}
if there is an action of $G$ on $M$, 
a $G$-equivariant map $\mu:M\to G$ (where $G$ acts on itself by
conjugation) and a $G$-invariant holomorphic two-form
$\omega\in \Omega^2(M)$ such that:

\noindent(QH1). 
$d\omega = \mu^*(\eta)$

\noindent(QH2).
For all $X\in \g$,
$\omega(v_X,\cdot\,) = \frac{1}{2}\mu^*(\Theta+\overline\Theta, X)
\in \Omega^1(M)$

\noindent(QH3).
At each point $m\in M$: 
$\ker \omega_m \cap \ker  d\mu = \{0\} \subset T_mM$. 

Here we have chosen a symmetric nondegenerate invariant bilinear form $(\ ,\ ):\g\otimes\g\to \IC$,
the Maurer--Cartan forms on $G$ are denoted $\Theta,\overline\Theta\in\Omega^1(G,\g)$ 
respectively (so in any representation 
$\Theta=g^{-1}dg, \overline\Theta=(dg)g^{-1}$), and
the canonical 
bi-invariant three-form on $G$ is
$\eta:= \frac{1}{6}([\Theta,\Theta],\Theta)$.
Moreover if $G$ acts on $M$, $v_X$ is the fundamental vector field of $X\in\g$; it is minus the tangent to the flow 
(so that the map $\g\to\Vect_M; X\mapsto v_X$ is a Lie algebra homomorphism). 

First we note that $\wh \cC$ is a complex manifold, in fact a smooth algebraic variety.
Let $\wt G$ denote the subvariety of $G\times \IP$ of pairs $(g,P)$ with $g\in P$ (this is the multiplicative Brieskorn--Grothendieck space if $P_0$ is a Borel).
There is a surjective map 
$$\pr: G\times P_0\to \wt G;\qquad (C,p)\mapsto (g,P)=(C^{-1}pC, C^{-1}P_0C)$$ 
whose fibres are precisely the orbits of a free action of $P_0$:
explicitly $q\in P_0$ acts 
on $G\times P_0$ as $q(C,p) = (qC, qpq^{-1})$.
Now choose a Levi decomposition $P_0 = LU$ of $P_0$ so that $U$ is the unipotent radical of $P_0$ and $L\cong P_0/U$. 
Consider the 
(locally closed) subvariety $\cC U\subset LU$ of $P_0$.
Since $\cC$ is a conjugacy class of $L$ and $P_0$ acts on $L$ via the projection $\pi : P_0\to L$, the conjugation action of $P_0$ on itself preserves $\cC U$.
Then $\pr$ restricts to $G\times \cC U$ and its image is $\wh \cC$, so that
$$\wh \cC \cong G\times_{P_0} \cC U$$
and we deduce $\wh \cC$ is a smooth complex algebraic variety. (Note that $\cC$ has a natural algebraic structure as a quotient of $L$, but will not be affine unless it is a semisimple conjugacy class).

Rather than prove directly that $\wh \cC$ is quasi-Hamiltonian we will use the well-known fact that $\cC$ is a quasi-Hamiltonian $L$-space and obtain $\wh \cC$ by reduction from a quasi-Hamiltonian $G\times L$ space, as follows.

Recall that $P_0$ acts freely on $G\times P_0$.
Let $\IM$ denote the quotient $G\times_U P_0 = (G\times P_0)/U$
by the subgroup $U\subset P_0$.
Thus $\IM$ has a residual  action of $L\cong P_0/U$ and also has a commuting action of $G$, from the action  
$g(C,p)= (Cg^{-1}, p)$ on $G\times P_0$.
Moreover there is a map $\mu:\IM\to G\times L$ induced from the $U$ invariant map
$$\wh \mu : G\times P_0 \to G\times L\qquad (C, p) \mapsto (C^{-1} p C, \pi(p)^{-1})$$
where $\pi : P_0\to L$ is the canonical projection.

\begin{thm} \label{thm: bigqh}
The space $\IM$ is a quasi-Hamiltonian $G\times L$ space with moment map $\mu$ and two-form $\omega$ determined by the condition 
$${\pr}^*(\omega) 
=\frac{1}{2} 
\left(\overline\gamma, \Ad_p\overline \gamma\right)
+\frac{1}{2} \left( \overline \gamma,\cp+\overline \cp\right)
\in \Omega^2(G\times P_0)
$$
where
$\pr$ is the projection $G\times P_0\to \IM$ and 
$\overline \gamma  = C^*(\overline\Theta)$,
$\cp  = p^*(\Theta),\bar\cp  = p^*(\bar\Theta)$.
\end{thm}

To deduce Theorem B from this we may perform the 
fusion $\IM\fus_L \cC$ with the conjugacy class $\cC\subset L$, and then perform the  quasi-Hamiltonian reduction by the free action of $L$.
The result $(\IM\fus_L \cC)\spq L$ may be identified immediately with $\wh \cC$.

In the special case $P_0=L=G$, the space $\IM$ is just the double $D(G)\cong G\times G$ of \cite{AMM}.
In general $\dim \IM = 2\dim P_0$ and the two-form $\omega$ on $\IM$ may be derived from the two-form $\omega_D$ on $D(G)$: one finds that the restriction of $\omega_D$ to  $G\times P_0$ (via the inclusion $P_0\subset G$) is basic for the $U$ action and descends to the two-form $\omega$ on $\IM$. That the result is again quasi-Hamiltonian requires proof of course.

\begin{rmk}
In the first instance the two-form $\omega$ was arrived at by actually computing what arose from the Hamiltonian loop group spaces related to resonant logarithmic connections on a disk, similarly to Section 4 of \cite{saqh}.
This computation led to the two-form from the double.
Note that in the case of $G=\GL_n(\IC)$ the quasi-Hamiltonian spaces $\wh \cC$ may be constructed differently, in terms of quivers (see \cite{yamakawa-mpa}), although even for $\GL_n$ the spaces $\IM$ do not seem to arise from quivers.
\end{rmk}

\begin{rmk}
Notice also that there are certain parallels with the Stokes phenomenon;
e.g. for the global moduli spaces, again one must fix a certain {\em union} of local gauge orbits to fix a symplectic leaf (in \cite{smid} this union
arose by fixing the {\em formal} gauge orbits). 
Also the spaces $\IM$ may be viewed as a tame analogue of the fission spaces of \cite{qf} (and again one may glue on more complicated spaces, not just conjugacy classes).
\end{rmk}

\pfms (of Theorem \ref{thm: bigqh}).
Write $\wh \omega = \pr^*(\omega)$. 
Let $U_-$ denote the unipotent radical of the parabolic opposite to $P_0$
(so that $U_-P_0$ is open in $G$).
(QH1) may be deduced directly from the result of \cite{qf} that $\gal:=G\times L \times U_-\times U$ is a quasi-Hamiltonian $G\times L$ space with moment map $\mu_\cA:(C, h, u_-, u)\mapsto (C^{-1} p C, h^{-1})\in G\times L$ where $p=u_-^{-1}hu$.
Then we can consider the embedding 
$$\iota:G\times P_0\hookrightarrow \gal
\qquad
(C, p)\mapsto(C, h, 1, u)$$
defined via the Levi decomposition $p=hu\in P_0$.
Thus $\wh \mu = \mu_\cA \circ \iota$ and moreover $\iota^*\Omega = \wh \omega$  where $\Omega$ is  the quasi-Hamiltonian two-form on $\gal$ from \cite{qf}. 
Then (QH1) follows immediately: 
$$\pr^*d\omega = d\pr^*\omega = d\iota^*\Omega = \iota^*d\Omega
= \iota^*\mu^*_\cA\eta = \wh\mu^*\eta = \pr^*\mu^* \eta$$
so that $d\omega=\mu^* \eta$ since $\pr$ is surjective on tangent vectors.
(QH2) is straightforward and left as an exercise.
(QH3) is trickier and we proceed as follows. %
 It is sufficient to show that at each point $m\in M:=G\times P_0$ the subspace $\ker \wh \omega \cap \ker d\wh \mu$ of the tangent space $T_mM$  is contained in the space of tangents to the $U$ action.
Thus choose $X\in T_mM$ and suppose that 
$X\in\ker(\wh \omega)\cap \ker(d\wh\mu)$.
Write $\wh\mu = (\mu_G, \mu_L)$ for the components of the moment map.
Since $X$ is in the kernel of $d\mu_L$ we have
$\ch'=0$ (here primes denote derivatives along $X$, so 
$\ch':=\langle h^*(\Theta_L),X\rangle$, where $\Theta_L$ is the Maurer-Cartan form on $L$).
Moreover $X$ being in the kernel of $d\mu_G$ amounts to the condition
$\bar\ga'+\cP'=p^{-1}\bar\ga' p$. Since $p=hu$ (and $\ch'=0$) this becomes
\beq\label{star}
\bar\ga' +\cU' = p^{-1}\bar\ga' p.
\eeq
(In general here the adjoint action of $g\in G$ on $X\in \g$ will be 
denoted $gXg^{-1}:=\Ad_gX$.)
Now we choose an arbitrary tangent vector $Y\in T_mM$ and denote 
derivatives along $Y$ by dots, 
so e.g.  $\dot\cP = \langle Y, \cP_m \rangle\in \Lie(P_0)$.
We then compute
$$2\wh\omega(X,Y) = 
2\left(
\bar\ga'
, \dot\cU\right)
+\left( 
u\bar\ga'u^{-1} + h^{-1}\bar\ga'h
, \dot\ch \right).
$$
\ppb{
\begin{align}\label{row1}
2\wh\omega(X,Y) &= 
\left(
p^{-1}\bar\ga'p 
- p(\bar\ga'+\cU')p^{-1}
 - \cU' , \dot{\bar\ga}\right)\\\label{row2}
&+\left(
\bar\ga' + p^{-1}\bar\ga'p
, \dot\cU\right)\\\label{row4}
&+\left( 
u\bar\ga'u^{-1} + h^{-1}\bar\ga'h
, \dot\ch \right).
\end{align} }
This should be zero for all $Y$; observe that each term on the right is really an independent condition on $X$.
From the first term we deduce the component of
$\bar\ga'$ in $\Lie(U_-)$ is zero.
The second term implies the $\Lie(L)$ component of $\bar\ga'$ is also zero.
Thus we find that $\bar\ga' \in \Lie(U)$, and we know $\ch'=0$ and  equation  \eqref{star} holds. But these three conditions characterise\footnote{To see this choose $X\in \Lie(U)$ consider the flow
$(C(t),p(t)) = \exp(Xt)\cdot (C, p)$. 
Thus (differentiating with respect to $t$) $\overline \ga' = C'C^{-1} = X\in \Lie(U)$ and similarly $\cP' = p^{-1}Xp - X$. But since $p=hu$ we see $h$ is constant and $\cU' = \cP'$, so \eqref{star} follows.} 
the tangents to the $U$ orbits on $G\times P_0$, so (QH3) follows.

\ppb{
Generally if $\A,\cB,\cC\in\Omega^1(M,\g)$ are $\g$-valued
holomorphic one-forms on a complex manifold $M$ then
$(\A,\cB)\in\Omega^2(M)$ and 
$[\A,\cB]\in\Omega^2(M,\g)$ are defined
by wedging the form parts and pairing/bracketing the Lie algebra parts
(so e.g. $(A\al,B\be)= (A,B)\al\wedge\be$ for
$A,B\in \g, \al,\be\in\Omega^1(M)$).
Define $\A\A := \frac{1}{2}[\A,\A]\in \Omega^2(M,\g)$ (which works out
correctly in
any representation of $G$  using matrix multiplication).
Then one has 
$d\Theta=-\Theta^2, d\overline\Theta = \overline\Theta^2$.
Define $(\A\cB\cC) = (\A, [\cB,\cC])/2\in\Omega^3(M)$ (which is
invariant under all permutations of $\A, \cB, \cC$).
}

\epfms

\begin{rmk}
In particular it follows that all the spaces $\wh \cC$ arise as certain moduli spaces of framed 
connections on a disc.
The precise statement is as follows.
Let $\Delta$ be a closed disc in the $z$ plane centered at zero. Replace $LG, \cA_\th, \wh \cP_\th$ by their analogues defined on all of $\Delta$ (rather than just germs at $0$). So e.g. now $LG = G(R)$ where $R$ is the ring of meromorphic functions of $\Delta$, having poles only at $0$. Choose a point $q$ on the boundary of $\Delta$ and let $\wh \cP^1_{\th}$ be the subgroup of $\wh \cP_\th$ of elements taking the value $1\in G$ at $q$.
Also let $\cA_\th(O)$ denote the subset of $\cA_\th$ of elements lying over a fixed orbit $O$
as in Theorem \ref{thm: main corresp} (and suppose $\cC, \phi$ are as defined there too).
Then 
$$\wh\cC \cong \cA_\th(O)/\wh\cP_\th^1$$
i.e. $\wh\cC$ is isomorphic to the space of connections on $\Delta$ lying over $O$ with a framing at $q$. Moreover the residual action of 
$\wh\cP_\th/\wh\cP_\th^1\cong G$ corresponds to the $G$ action on $\wh \cC$.
\end{rmk}

\section{Cleaner statement}\label{sn: buildings}

A cleaner Riemann--Hilbert statement arises if we also allow the weight $\th$ to vary in the correspondence, but for this we 
need to define the notion of a weighted parahoric subgroup, analogous to the notion of weighted parabolic subgroup.
This leads directly to the definition of  Bruhat--Tits building.

First define the partly extended affine Weyl group to be $\wh W = N(\cK)/T(\cO)\cong W\sdp X_*(T)$ where $X_*(T)$ is the cocharacter lattice, which we think of either as the set of 1 parameter subgroups of $T$, or as the kernel of $\exp(2\pi i \,\cdot):\lt\to T$ (an element $\la$ of this kernel corresponds to the one parameter subgroup 
$\varphi = z^\la$).
Here $N\subset G$ is the normaliser of $T$ in $G$ and 
$W = N/T$ is the finite Weyl group, which acts naturally on $\lt_\IR$ (via the adjoint action of $N$).
Note that $X_*(T)\cong T(\cK)/T(\cO)$ and by convention $X_*(T)$  acts on $\lt_\IR$ via $z^\la \cdot \theta = \theta -\la$ 
(this is a standard convention, but beware it agrees with our conventions concerning gauge transformations only if we identify $\th$ with minus the residue of the connection $-\th dz/z$.)
These two actions combine to give an action of $\wh W$ on $\lt_\IR$.

\begin{defn}
A {\em weighted parahoric subgroup} of $LG$ is an equivalence class of elements $(g,\th)\in LG \times \lt_\IR$ where 
$$(g,\th)\sim (g',\th')$$
if $\th' =w \th$ for some $w\in \wh W$ and 
$g^{-1}g'\wh w \in \wh \cP_\th$
for some lift $\wh w$ of $w$ to $N(\cK)\subset LG$.
\end{defn}

\ppb{
Here the finite Weyl group $W=N(T)/T$ acts on $\lt_\IR$ in the usual way and a one-parameter subgroup $\varphi=z^\la\in X_*(T)$ acts by mapping $\th\in \lt_\IR$ to $\th-\la$. (This is a standard convention, but beware it agrees with our conventions concerning gauge transformations only if we identify $\th$ with the connection $-\th dz/z$.)
}
This is the standard definition of the (extended) Bruhat--Tits building $\cB(LG)= (LG\times \lt_\IR)/\sim$ of $LG$ \cite{BrTits-I} p.170.
Thus we are saying a weighted parahoric is a {\em point} of the building. 
(It seems one usually views the building as a simplicial complex and rarely regards its points in this sense.)
Note that $LG$ acts naturally on $\cB(LG)$ via left multiplication on $LG$.

\begin{lem}
A weighted parahoric $p\in \cB(LG)$ canonically determines a parahoric subgroup $\cP_p\subset LG$ and a space of connections $\cA_p\subset \cA$.
\end{lem}
\pf
Suppose $p$ is in the equivalence class of $(g,\th)\in LG\times \lt_\IR$, and  $(g',\th')\sim (g,\th)$ with
$\th' = w\th$. Choose a lift $\wh w$ of $w$ to $N(\cK)=N(\IC)\sdp T(\cK)\subset LG$.
We may check directly that 
$\cP_{\th'} = \wh w\cP_\th \wh w^{-1}$ and that $\cA_{\th'} = \wh w[\cA_\th]$.
The first claim then follows since 
$\wh\cP_\th$ normalises $\cP_\th$: $\cP_p := g\cP_\th g^{-1}$ is well defined.
Secondly we should check that $\cA_p := g[\cA_\th]$ depends only on the equivalence class of $p$.
But by Lemma \ref{lem: wP preserves Ath}  
$g^{-1}g'\wh w$ 
preserves $\cA_\th$ so $\cA_p = g'[\cA_{\th '}]$.
\epf

\begin{rmk}
Note there is an embedding 
$\lt_\IR\hookrightarrow \cB(LG); \th\mapsto [(1,\th)]$ (whose image is the standard apartment) and one may then confirm (see Lemma \ref{lem: parah stabs}) that $\wh \cP_\th$ is exactly the stabiliser in $LG$ of $\th\in \cB(LG)$.
It follows in general that $\cP_p$ is the identity component of $\Stab_{LG}(p)$.
\end{rmk}

Thus it makes sense to consider pairs $(A, p)$ where $p\in \cB(LG)$ is a weighted parahoric and $A\in \cA_p$ is a compatible connection.
It follows from the lemma  that the  loop group $LG$ acts on the set of such pairs: $g(A,p) = (g[A], g(p))$.

The corresponding  monodromy data consists of pairs 
$(M, b)\in G\times \IB(G)$ with $M\in P_b$.
Here $\IB(G)$ is the space of weighted parabolic subgroups of $G$. 
A point of $\IB(G)$ consists of a parabolic $P\subset G$ and a set of weights for $P$ (as defined earlier). This can be rephrased 
to parallel the definition of $\cB(LG)$ as follows.

\begin{defn}
A {\em weighted parabolic subgroup} of $G$ is an equivalence class of elements $(g,\th)\in G \times \lt_\IR$ where 
$(g,\th)\sim (g',\th')$
if $\th' =w \th$ for some $w\in W$ in the Weyl group and 
$g P_\th g^{-1} =  g' P_{\th'} (g')^{-1}\subset G$ 
(i.e. $g^{-1}g'\wh w\in P_\th$ for some lift $\wh w\in N(\IC)$ of $w$).
\end{defn}

Thus we can define $\IB(G) =  (G \times \lt_\IR)/\!\sim$ and note that 
$b\in \IB(G)$ determines a parabolic subgroup $P_b = gP_\th g^{-1}\subset G$. (Beware this is not the spherical building of $G$, it is more like the cone over the spherical building; if we choose a maximal compact subgroup $K\subset G$ then one may identify 
$\IB(G)\cong i\Lie(K)\subset \g$.)\footnote{
Similarly it seems one may identify $\cB(LG)$ with a space of $K$-connections, although we will not use this viewpoint.}
In any case, basically as a corollary of Theorem  \ref{thm: main corresp} we find:

\begin{cor}
There is a canonical bijection between $LG$ orbits of tame parahoric connections and $G$ orbits of enriched monodromy data:
$$\Bigl\{ (A,p) \st p\in \cB(LG), A\in \cA_p\Bigr\}/LG \quad \cong\quad 
 \Bigl\{ (M,b) \st b\in \IB(G), M\in P_b\Bigr\}/G.$$
\end{cor}
\pf
Given $(A,p)$ we may act by $LG$ to move $p$ to a point $\th$ of the standard apartment, and thus suppose $A\in \cA_\th$ and $p=\th\in \lt_\IR$. 
We may further assume $A$ is in normal form.
Then we may obtain data $M, \phi$ as usual, with $M\in P_\phi$, i.e. a point of the right-hand side, with $b=\phi$.
We should check that the $G$-orbit of $(M,b)$ only depends on the $LG$ orbit of $(A,p)$: firstly this is clear if we only move $(A,p)$ by an element of $\wh \cP_\th$ (so $p=\th$ does not move) by Lemma 
\ref{lem: orbit idpce} 
and Theorem \ref{thm: main corresp}.
Secondly we should examine what happens if we act by an element $g$ of $N(\cK)= N(\IC)\sdp T(\cK)$
(since any other element of the $LG$ orbit of $(A,p)$ above the standard apartment will arise in this way).
We may write $g =  h z^\mu t$ with 
$h\in N(\IC), \mu\in X_*(T), t\in T(\cO)$.
Since $t\in \wh \cP_\th$ we may assume $t=1$ here.
Set $A' = g[A], \th' = g\cdot \th = \Ad_h(\th-\mu)$.
It is straightforward to check that $A'$ is again in normal form: indeed suppose $A = (\tau+\si +\sum A_{i\al}z^i) dz/z$ with $\al\in \cR\cup\{0\}$ and $A_{i\al}\in \g_\al$ (or in $\lt$ if $\al=0$),  then 
$$hz^\mu[A] = \Ad_h(\tau+\mu + \si +\sum A_{i\al}z^{i+\al(\mu)})\frac{dz}{z}.$$ 
The key point then is that $\tau' = \Ad_h(\tau+\mu)$ so that
$$\phi'  = \tau'+\th' = \Ad_h(\tau+\mu+\th-\mu)=\Ad_h(\phi)$$
so that $\phi$ only moves via the {\em finite} Weyl group $W$.
The corresponding fundamental solutions are of the form $z^\tau z^R$ and $z^{\tau'}z^{R'} = hz^{\tau+\mu} z^R h^{-1}$ so it is clear that the monodromies etc. are related by the action of $h$.
This shows the map from left to right is well-defined. Surjectivity follows from Theorem \ref{thm: main corresp}. Injectivity also largely follows from Theorem \ref{thm: main corresp}, but it remains to check that orbits with inequivalent $\th$ map to different points.
But this follows from that fact that $(M,b)$ determines the $\wh W$ orbit of $\th\in\lt_\IR$---indeed suppose 
we act by $G$ so that 
$M\in P_\phi$, with $\phi\in \lt_\IR$ determined upto the action of $W$. Then let $d\in T$ be any element conjugate to the semisimple part of $\pi(M)\in L=C_G(\phi)$, so that $d= \exp(2\pi i (\tau+\si))$ with $\tau$ determined up to the addition of an element of $X_*(T)$.
This yields one choice of $\th = \phi-\tau$ and the others are determined by making different choices---i.e. via the action of $\wh W$.
\epf

\section{Other directions.}

\ppb{First it looks to be a formality (although we have not yet checked the details) to deduce a global Riemann--Hilbert correspondence over smooth complex  algebraic curves $\Si$, %
involving connections on torsors for parahoric (Bruhat--Tits) group schemes $\G\to\Si$, such that locally 
$\G$ looks like a parahoric subgroup $\cP$ of the local loop group
and at all but finitely many points %
 $\G$ looks like $G(\cO)$.
Such torsors are studied (in more generality, but not with connections or weights) in \cite{pap-rap}.
Also
it seems it will be possible to extend the nonabelian Hodge correspondence to the present context (i.e. the correspondence between such connections and Higgs bundles);
The correspondence of the parameters between connections and Higgs bundles will be as in Simpson's table \cite{Sim-hboncc} p.720---basically the parameters are rotated, and this now generalises directly.
}

First it looks to be possible to extend the nonabelian Hodge correspondence to the present context (i.e. the correspondence on a smooth algebraic curve $\Si$ between such connections and Higgs bundles, under stability conditions);
The correspondence of the parameters will be as in Simpson's table \cite{Sim-hboncc} p.720---basically the parameters are rotated, and this now generalises directly.
In our notation this table is:

  \begin{center}
  \begin{tabular}{| c || c | c | c |}
    \hline
     & \text{Dolbeault} & \text{De\! Rham} & \text{Betti}\\ \hline
    \text{weights }$\in \lt_\IR$ & 
$-\tau$ & $\th$ & $\phi=\tau+\th$ \\ \hline
    \text{``eigenvalues''}$\in \lt_\IC, \lt_\IC, \text{T}(\IC)$ & 
$-\frac{1}{2}(\phi+\si)$ & $-(\tau+\si)$ & $\exp(2\pi i (\tau+\si))$ \\
    \hline
  \end{tabular}
  \end{center}
where the columns correspond to Higgs bundles, connections and monodromy data respectively.\footnote{Beware that we use the opposite conventions to \cite{Sim-hboncc} for connections on vector bundles 
($d-A$ rather than $d+A$ in local trivialisations)---this explains the sign in the middle of the bottom row, that does not appear elsewhere in the present article.} Observe for example that the eigenvalues of the Higgs field will only vary under the finite Weyl group, as expected.
This global correspondence is probably best phrased in terms of torsors for parahoric (Bruhat--Tits) group schemes $\G\to\Si$, such that locally 
$\G$ looks like a parahoric subgroup $\cP$ of the local loop group
and at all but finitely many points %
 $\G$ looks like $G(\cO)$.
Such torsors have been  studied recently (in more generality, but not with connections or weights) in \cite{pap-rap}.
On the other hand quasi-parahoric Higgs bundles (i.e. without the weights)
have been studied algebraically recently by Yun \cite{yun2},
and, in effect, the local picture of such Higgs bundles was studied by Kazhdan--Lusztig \cite{kazlus88} in 1988. (Corollary E is related to the De\! Rham and Betti analogues of this.) 
It might also be profitable (in the case of rational weights) to relate  the parahoric viewpoint here to the ``ramified'' approach of Balaji et al \cite{BBN03} (see also Seshadri \cite{sesh08}), although they have not  considered the analogue of logarithmic connections it seems.

\appendix

\section{Extra proofs}

\begin{lem}\label{lem: parah stabs}
For any $\th\in \lt_\IR$, the group $\wh \cP_\th$
is the stabiliser in $LG$ of $p=[(1,\th)]\in \cB(LG)$.
\end{lem}
\pf
Clearly $\wh \cP_\th$ does stabilise $p$.
Conversely if $g\in LG$ stabilises $p$ then
$(g,\th)\sim (1,\th)$ so that $g^{-1}\wh w \in \wh\cP_\th$ for some $\wh w\in N(\cK)$ such that $w(\th)= \th$
(where $w$ is the image of $\wh w$ in $\wh W$).
Thus it is sufficient to show that all such elements $\wh w$ are in $\wh \cP_\th$.
Thus we should check that $z^\th \wh w z^{-\th}$ has a limit as $z\to 0$ along any ray.
We may write $\wh w = hz^\la t$ with $h\in N(\IC), \la\in X_*(T), t\in T(\cO)$.
Thus 
$$z^\th \wh w z^{-\th} = z^\th h z^{-\th+\la}t 
= hz^{\Ad^{-1}_h(\th)}    z^{-\th+\la}t. $$
But the condition that $w\th = \th$ means 
$\Ad_h(\th-\la) = \th$, so that $\Ad^{-1}_h(\th) = \th-\la$, and the above expression reduces to $ht$, which clearly has a limit as $z\to 0$.
\epf

\renewcommand{\baselinestretch}{1}              %
\normalsize
\bibliographystyle{amsplain}    \label{biby}
\bibliography{../thesis/syr} 

\vspace{0.5cm}   
\'Ecole Normale Sup\'erieure et CNRS, 
45 rue d'Ulm, 
75005 Paris, 
 France

www.math.ens.fr/$\sim$boalch

boalch@dma.ens.fr \qquad \qquad \qquad \qquad \qquad \quad\,\,  

\end{document}

%% file: macros.tex
\usepackage{verbatim}  

\newenvironment{ppb}[1]
{\ \!\!\!\!\!\!\!\!\!\!\!\!\!\!\!\!\!\!\!\!\!\!\!\!\!\!\!\!\!\!\!\!\!\!\!\!\!\!\!\! {\bf PPB------------------------------------------------------------------------------------------------PPB}\newline \tiny {#1}
\  \newline\normalsize\phantom{f}\!\!\!\!\!\!\!\!\!\!\!\!\!\!\!\!\!\!\!\!\!\!\!\!\!\!\!\!\!\!\!\!\!\!\!\!\!\!\!\! {\bf PPB------------------------------------------------------------------------------------------------PPB}\newline}{}

\newenvironment{ssn}[1]{{\ \newline \large\bf\hspace{-\parindent}{#1.}}}{}

\newcommand\beq{\begin{equation}}
\newcommand\eeq{\end{equation}}
\newcommand\bal{\begin{align*}}
\newcommand\eal{\end{align*}}   
\newcommand\bmx{\left(\begin{matrix}}
\newcommand\emx{\end{matrix}\right)}
\newcommand\bsmx{\left(\begin{smallmatrix}}
\newcommand\esmx{\end{smallmatrix}\right)}

\newcommand{\spq}{/\!\!/}

\newcommand{\st}{\ \bigl\vert\ }

\providecommand{\Rad}{\text{\rm Rad}}

\def\part#1{\frac{\partial\phantom{q}}{\partial#1}}

\newcommand {\flb}{\lbrack\!\lbrack}
\newcommand {\frb}{\rbrack\!\rbrack}

\newcommand{\sdp}{{\ltimes}}

\newcommand{\fus}{\circledast}

\newcommand{\Lie}{{\mathop{\rm Lie}}}

\newcommand{\Vect}{{\mathop{\rm Vect}}}             

\newcommand{\gal}{{\ _G\cA_L}}

\newcommand{\Ad}{{\mathop{\rm Ad}}}
\newcommand{\ad}{{\mathop{\rm ad}}}

\DeclareMathOperator{\pr}{pr}

\newcommand{\SL}{{\mathop{\rm SL}}}

\newcommand{\GL}{{\mathop{\rm GL}}}
\newcommand{\PGL}{{\mathop{\rm PGL}}}

\renewcommand{\ker}{\mathop{\rm Ker}}

\newcommand{\End}{\mathop{\rm End}}

\newcommand{\Stab}{{\mathop{\rm Stab}}}

\newcommand{\IB}{\mathbb{B}}
\newcommand{\IC}{\mathbb{C}}
\newcommand{\ID}{\mathbb{D}}

\newcommand{\IL}{\mathbb{L}}
\newcommand{\IM}{\mathbb{M}}

\newcommand{\IP}{\mathbb{P}}                                     
                           
\newcommand{\IR}{\mathbb{R}}

\newcommand{\IZ}{\mathbb{Z}}

\newcommand{\A}{\mathcal{A}}
\newcommand{\cA}{\mathcal{A}}
\newcommand{\cB}{\mathcal{B}}
\newcommand{\cC}{\mathcal{C}}
\newcommand{\cD}{\mathcal{D}}

\newcommand{\G}{\mathcal{G}}

\newcommand{\ch}{\hslash}    
\newcommand{\cI}{\mathcal{I}}

\newcommand{\cK}{\mathcal{K}}
\newcommand{\cL}{\mathcal{L}}

\newcommand{\cO}{\mathcal{O}}
\newcommand{\orbit}{O}
\newcommand{\cP}{\mathcal{P}}
\newcommand{\cp}{\mathcal{P}}

\newcommand{\cR}{\mathcal{R}}

\newcommand{\cU}{\mathcal{U}}

\newcommand{\g}{       \mathfrak{g}     }

\newcommand{\lb}{       \mathfrak{b}            }

\newcommand{\n}{\mathfrak{n}}
\newcommand{\lt}{\mathfrak{t}}
\newcommand{\lh}{\mathfrak{h}}

\renewcommand{\sl}{       \mathfrak{sl}     } 

\newcommand{\h}{\mathfrak{h}}

\newcommand{\lp}{\mathfrak{p}}
\newcommand{\llp}{\wp}
\newcommand{\lln}{\mathfrak{ln}}
\renewcommand{\ll}{\mathfrak{l}}

\newcommand{\wt}{\widetilde}

\newcommand{\wh}{\widehat}

\newcommand{\al}{\alpha}

\newcommand{\be}{\beta}
\newcommand{\ga}{\gamma}

\newcommand{\la}{\lambda}

\newcommand{\si}{\sigma}

\newcommand{\Si}{\Sigma}
\renewcommand{\th}{\theta}

\renewcommand{\bar}{\overline}

\makeatletter
 \newlength{\typesize}
\makeatother

\newlength{\vvoff}
\newlength{\hhoff}

\newcommand{\pf}{\begin{bpf}}

\newcommand{\pfms}{\begin{bpfms}}
\newcommand{\epf}{\end{bpf}\hfill$\square$\\}           
\newcommand{\epfms}{\end{bpfms}\hfill$\square$\\}       

\newcommand{\idea}{\begin{bidea}}

\newcommand{\eidea}{\end{bidea}\hfill$\square$\\}           

\newcommand{\sk}{\begin{bsk}}    

\newcommand{\esk}{\end{bsk}\hfill$\square$\\}           
\newcommand{\sketch}{\begin{bsketch}}

\newcommand{\esketch}{\end{bsketch}\hfill$\square$\\}



%% file: macros-thm1.tex
\newtheorem {hypo}{\bf\hspace{-\parindent}Hypothesis}
\newtheorem {thm}[hypo]{Theorem}   

\newtheorem {cor}[hypo]{Corollary}
\newtheorem {lem}[hypo]{Lemma}

\newtheorem {defn}[hypo]{Definition}

\theoremstyle{remark}\newtheorem{rmk}[hypo]{Remark}